\documentclass[10pt]{article}
\usepackage{pb-diagram}
\usepackage{amsmath, amsfonts}
\usepackage{amssymb}
\usepackage{amscd}
\usepackage{hyperref}
\newtheorem{Df}{Definition}[section]
\newtheorem{Te}[Df]{Theorem}
\newtheorem{Po}[Df]{Proposition}
\newtheorem{Cr}[Df]{Corollary}
\newtheorem{Lm}[Df]{Lemma}
\newtheorem{Ca}[Df]{Claim}
\newtheorem{Cn}[Df]{Conjecture}
\newtheorem{Ex}[Df]{Example}
\newtheorem{Rm}[Df]{Remark}

\newcommand{\Bdf}{\begin{Df}}
\newcommand{\Edf}{\end{Df}}
\newcommand{\Bte}{\begin{Te}}
\newcommand{\Ete}{\end{Te}}
\newcommand{\Bpo}{\begin{Po}}
\newcommand{\Epo}{\end{Po}}
\newcommand{\Bcr}{\begin{Cr}}
\newcommand{\Ecr}{\end{Cr}}
\newcommand{\Blm}{\begin{Lm}}
\newcommand{\Elm}{\end{Lm}}
\newcommand{\Bca}{\begin{Ca}}
\newcommand{\Eca}{\end{Ca}}
\newcommand{\Bcn}{\begin{Cn}}
\newcommand{\Ecn}{\end{Cn}}
\newcommand{\Bex}{\begin{Ex}}
\newcommand{\Eex}{\end{Ex}}
\newcommand{\Brm}{\begin{Rm}}
\newcommand{\Erm}{\end{Rm}}
\newcommand{\Bdm}{{\it Proof.}\ }
\newcommand{\Edm}{\rule{2mm}{2mm}}

\newcommand{\bpf}{\begin{proof}}
\newcommand{\epf}{\end{proof}}

\newcommand{\ff}{\ensuremath{\mathbb{F}}}

\newcommand{\zz}{\ensuremath{\mathbb{Z}}}
\newcommand{\kk}{\ensuremath{\mathbb{K}}}

\newcommand{\pgq}{\geqslant}

\newcommand{\mo}{\ensuremath{\mathfrak{s}}}

\newcommand{\me}{\ensuremath{\mathfrak{e}}}

\newcommand{\id}{\operatorname{id}\nolimits}

\newcommand{\car}{\operatorname{char}\nolimits}

\newcommand{\kkmod}{\ensuremath{\kk}\text{-\textsf{\upshape Mod}}}
\newcommand{\fmod}{\ensuremath{\ff}\text{-\textsf{\upshape Mod}}}
\newcommand{\almod}{\ensuremath{A}\text{-\textsf{\upshape Mod}}}
\newcommand{\aemod}{\ensuremath{A^e}\text{-\textsf{\upshape Mod}}}

\newcommand{\blmod}{\ensuremath{B}\text{-\textsf{\upshape Mod}}}

\newcommand{\abimod}{\ensuremath{A}\text{-\textsf{\upshape Bimod}}}

\newcommand{\Hom}{\operatorname{Hom}\nolimits}
\newcommand{\Und}{\operatorname{Und}\nolimits}

\newcommand{\RHom}{\mathop{\rm RHom}\nolimits}

\newcommand{\HH}{\operatorname{HH}\nolimits}
\newcommand{\HK}{\operatorname{HK}\nolimits}
\newcommand{\cupk}{\underset{K}{\smile}}
\newcommand{\capk}{\underset{K}{\frown}}

\title{\bf{Calabi-Yau property in derived Koszul calculus}}
\author{Roland Berger\footnote{\href{roland.berger@univ-st-etienne.fr}{roland.berger@univ-st-etienne.fr}}
\and Jun Maillard \footnote{\href{jun.maillard@matfyz.cuni.cz}{jun.maillard@matfyz.cuni.cz}}}
\date{}

\begin{document}

\maketitle

\begin{abstract}
A Poincar\'e Van den Bergh duality theorem for strong Kc-Calabi-Yau algebras was obtained by R. Taillefer and the first author~\cite[Theorem 5.11]{bt:kcpa} under the assumption that the derived functors of functors involved in the statement exist. We prove the existence of these derived functors by showing that the dg category defining the derived Koszul calculus is isomorphic to a dg category of dg modules over a dg algebra. Therefore we get a definition of strong Kc-Calabi-Yau algebras and a corresponding duality theorem without any existence assumption. We prove that a polynomial algebra is strong Kc-Calabi-Yau.
\end{abstract} 

\section{Introduction} \label{intro}

Derived categories are systematically used as a fundamental framework in several important subjects. An emblematic case is given by the derived algebraic geometry, providing new approaches in algebraic geometry with a wide range of applications~\cite{toen:dag}. Cluster algebras and representation theory of quivers constitute another promising subject~\cite{keller:clustdercat}, including Ginzburg dg algebras developed from the seminal paper of V. Ginzburg on Calabi-Yau algebras~\cite{vg:cy}.

As a common feature of derived categories appearing in various contexts, the basic objects are dg algebras $\mathcal{A}$, dg $\mathcal{A}$-modules and abelian categories $\mathcal{M}$, giving rise to derived categories $\mathcal{D}(\mathcal{A},\mathcal{M})$. In our paper, we follow the recent book on derived categories by A. Yekutieli~\cite{yeku:dercat}. This book offers a self-contained text on the fundamentals of the theory, leading in great detail to a definition of derived categories $\mathcal{D}(\mathcal{A},\mathcal{M})$ and to existence theorems for derived functors.

Yelutieli's book contains an application to dualizing complexes -- a classical topic in algebraic geometry -- suitably extended to noncommutative dualizing complexes over noncommutative dg rings, following an approach by M. Van den Bergh~\cite{vdb:existence} developed by A. Yekutieli and J. J. Zhang (see~\cite{yeku:dercat} for references). In this way, derived categories appear as a natural bridge between commutative algebra and noncommutative algebra.

A standard theorem asserts that the algebra of multivector fields on a manifold $M$ is isomorphic to the Hochschild cohomology algebra of the algebra $A=C^{\infty}(M)$. The scope of this theorem was widely enlarged by M. Kontsevich in his formality theorem, see e.g.~\cite{lpv:poisson}. Here there is a bridge between differential geometry and homological algebra, which permits to define by analogy a differential calculus on noncommutative associative algebras, called Tamarkin-Tsygan calculus~\cite{dgt:variant, tt:calculus, tl:bvcy}, also named calculus structure~\cite{eu:calculus, eusched:cyfrob}.  

Koszul calculus was introduced by T. Lambre, A. Solotar and the first author in 2018~\cite{bls:kocal} in order to examine what happens for Tamarkin-Tsygan calculus when the bar resolution is replaced by the Koszul complex. In 2020, Koszul calculus was extended to quiver algebras with homogeneous quadratic relations (called \emph{quadratic quiver algebras}) by R. Taillefer and the first author~\cite{bt:kcpa}. This generalization, applied to preprojective algebras, reveals some interesting interactions with the Poincar\'e Van den Bergh duality~\cite{vdb:dual}. From the case of preprojective algebras, a new Calabi-Yau property was defined for quadratic quiver algebras in~\cite{bt:kcpa}, called \emph{Koszul complex Calabi-Yau property}, abbreviated in \emph{Kc-Calabi-Yau property}.

Explicitly, a quadratic quiver algebra $A$ is said to be $n$-Kc-Calabi-Yau if the $A$-bimodule Koszul complex $K(A)$ of $A$ has finite length $n$ and if 
$$\RHom_{A^e}(K(A), A^e) \cong K(A)[-n]$$
in the bounded derived category $\mathcal{D}^b(\abimod)$ of $A$-bimodules. As usual, $A^e:=A\otimes_{\ff} A^{op}$ where $\ff$ is the base field.

If $A$ is $n$-Kc-Calabi-Yau and Koszul, then $A$ is homologically smooth and $K(A) \cong A$ in $\mathcal{D}^b(\abimod)$, so that one has
$$\RHom_{A^e}(A, A^e) \cong A[-n] \ \mathrm{in} \ \mathcal{D}^b(\abimod),$$
meaning that $A$ is $n$-Calabi-Yau in Ginzburg's sense~\cite{vg:cy}. In~\cite{bt:kcpa}, it is proved that the preprojective algebra of a connected graph distinct from A${_1}$ and A${_2}$ is 2-Kc-Calabi-Yau. Such a preprojective algebra is known to be not 2-Calabi-Yau in Ginzburg's sense when the graph is Dynkin ADE, because the minimal projective resolution of $A$ has then infinite length so that $A$ is not homologically smooth. If the graph is not Dynkin ADE, the preprojective algebra is Koszul.

In Koszul calculus, the Hochschild cochain dg algebra $\Hom_{A^e}(B(A), A)$ of the Tamarkin-Tsygan calculus, where $B(A)$ is the bar resolution of $A$, is replaced by \emph{the Koszul cochain dg algebra of $A$}, denoted by $\tilde{A}$. By definition,
\begin{equation} \label{tildea1}
\tilde{A}=\Hom_{A^e}(K(A), A)
\end{equation}
is the space of Koszul cochains of $A$ with coefficients in $A$. This is a dg algebra whose grading is given by the cohomological degree of Koszul cochains, differential is the Koszul differential $b_K$ and multiplication is the Koszul cup product $\underset{K}{\smile}$. The elements of $\tilde{A}$ should be thought of as differentials on $A$. If $A$ is a polynomial algebra in $n$ variables, the algebra $\tilde{A}$ is isomorphic to the exterior algebra $\Omega_A^{\ast}$ constructed on the K\"{a}hler differentials of $A$~\cite{weib:homo}, as we shall see it in Theorem \ref{dualityiso}.

A fundamental fact is that the Koszul complex $K(A)$ is a dg $\tilde{A}$-bimodule compatible with its structure of $A$-bimodule~\cite[Proposition 3.6]{bt:kcpa}. Following the Yekutieli terminology~\cite{yeku:dercat} applied to dg bimodules instead of left dg modules, we say that $K(A)$ is a dg $\tilde{A}$-bimodule in the abelian category $\abimod$ of $A$-bimodules.

In the same spirit, notations used in the Yekutieli formalism for left dg modules will be adapted in our paper for dg bimodules. For example, we denote by $\mathcal{C}(\tilde{A},\abimod)$ the dg category of dg $\tilde{A}$-bimodules in $\abimod$ and by $\mathcal{D}(\tilde{A},\abimod)$ its derived category. The strict category of $\mathcal{C}(\tilde{A},\abimod)$ is $\mathcal{C}_{str}(\tilde{A},\abimod)$ while the homotopy category is $\mathcal{K}(\tilde{A},\abimod)$, giving rise to canonical functors
$$P: \mathcal{C}_{str}(\tilde{A},\abimod) \longrightarrow \mathcal{K}(\tilde{A},\abimod),$$
$$Q: \mathcal{K}(\tilde{A},\abimod) \longrightarrow \mathcal{D}(\tilde{A},\abimod).$$
Definitions and basic results on derived categories~\cite{yeku:dercat} are presented in Section \ref{generalresult} below. Let us introduce the following definition.

\Bdf \label{derkoscal1}
Let $A$ be a quadratic quiver algebra. The data of the categories $\mathcal{C}(\tilde{A},\abimod)$ and $\mathcal{D}(\tilde{A},\abimod)$, and of the functor
$$\tilde{Q}=Q \circ P : \mathcal{C}_{str}(\tilde{A},\abimod) \longrightarrow \mathcal{D}(\tilde{A},\abimod)$$
is called the derived Koszul calculus of $A$.
\Edf

Our first result concerning on derived Koszul calculus is the following proposition. This result was not observed in~\cite{bt:kcpa}. It is essential for defining strong Kc-Calabi-Yau algebras and for expressing a Poincar\'e Van den Bergh duality theorem as a Koszul cap product by a fundamental class, \emph{without any existence assumption of the involved derived functors}. 

\Bpo \label{fundaprop1}
Let $A$ be a quadratic quiver algebra. The category $\mathcal{C}(\tilde{A},\abimod)$ has enough $K$-injectives and enough $K$-projectives.
\Epo

This proposition, stated as Proposition \ref{fundaresult} in Section \ref{derkoscalA} below, will be a consequence of Theorem \ref{mixedsitutheo}, a new general result in the Yekutieli formalism whose statement is the following.

\Bte \label{anewgenresult}
Assume that $\kk$ is a nonzero commutative base ring, $A$ is a central $\kk$-ring and $B$ is a dg central $\kk$-ring. The category $\mathcal{C}(B, \almod)$ has enough $K$-injectives and enough $K$-projectives.
\Ete

Section \ref{generalresult} is devoted to prove this theorem. Therefore, as explained in Section \ref{derkoscalA}, the derived functor
\begin{equation} \label{bendohomderived}
\RHom_{A^e}(-, A^e):\mathcal{D}^b(\tilde{A},\abimod)^{op} \longrightarrow \mathcal{D}^b(\tilde{A},\abimod)
\end{equation}
exists, where the upperscript $b$ corresponds to bounded complexes.

At the end of Section \ref{derkoscalA}, we shall be ready to reformulate the definition of strong Kc-Calabi-Yau algebras \emph{without any existence assumption}. Our motivation for defining strong Kc-Calabi-Yau algebras is that the preprojective algebra of a connected graph distinct from A${_1}$ and A${_2}$ appears naturally as strong 2-Kc-Calabi-Yau~\cite[Proposition 4.7]{bt:kcpa} including an existence assumption, now unnecessary.

\Bdf \label{strongCY}
Let $A$ be a quadratic quiver algebra. We say that $A$ is strong $n$-Kc-Calabi-Yau if the complex $K(A)$ has length $n$ and if
$$\RHom_{A^e}(K(A), A^e) \cong K(A)[-n]$$
in the bounded derived category $\mathcal{D}^b(\tilde{A},\abimod)$.
\Edf

Then $A$ is $n$-Kc-Calabi-Yau by forgetting the actions of $\tilde{A}$. Concluding Section \ref{derkoscalA}, we shall obtain a Poincar\'e Van den Bergh duality theorem for strong Kc-Calabi-Yau algebras \emph{without any existence assumption}. It is a duality between Koszul cohomology and Koszul homology expressed as a cap product by a fundamental class, inducing a duality between their higher versions indicated by the subscript (upperscript) $hi$. 

\Bte \label{strongvdbduality1}
Let $A$ be a quadratic quiver algebra. Assume that $A$ is strong $n$-Kc-Calabi-Yau. Then there is a class $c$ in $\HK_n(A)$, called fundamental class, such that  
$$c \underset{K}{\frown} -  : \HK^{\bullet}(A) \rightarrow \HK_{n-\bullet}(A)$$
is an isomorphism of $\HK^{\bullet}(A)$-bimodules, inducing an isomorphism of $\HK^{\bullet}_{hi}(A)$-bimodules from $\HK^{\bullet}_{hi}(A)$ to $\HK^{hi}_{n-\bullet}(A)$. For all $\alpha \in \HK^p(A)$, one has $c \underset{K}{\frown} \alpha = (-1)^{np} \alpha \underset{K}{\frown}c$.
\Ete

Section \ref{polynomials} is devoted to prove that a polynomial algebra in $n$ variables is strong $n$-Kc-Calabi-Yau. It is based on a duality isomorphism stated below. This isomorphism is reminiscent of the contraction map (inner product) used in differential geometry~\cite{lpv:poisson}. Here the quiver has a single vertex and $n$ arrows, so that $k=\ff$. Fix a basis $\{x_1, \ldots, x_n\}$ of the $k$-vector space $V$ and define a fundamental class $c$ of $A$ by $c=1_A \otimes (x_1 \wedge \ldots \wedge x_n)$ in $A\otimes W_n$. 

\Bte \label{dualityisomorphism}
Let $V$ be a $k$-vector space of finite dimension $n$ and $A=S(V)$ be the symmetric algebra of $V$. Let $M$ be an $A$-bimodule. For each Koszul $p$-cochain $f$ with coefficients in $M$, we define the Koszul $(n-p)$-chain $\theta_M (f)$ with coefficients in $M$ by
\begin{equation} \label{definitiontheta}
\theta_M (f) = c \underset{K}{\frown} f.
\end{equation}
Then the linear map $\theta_M: \Hom(W_{\bullet},M) \rightarrow M\otimes W_{n-\bullet}$ is a $(-n)$-degree complex isomorphism of dg $\tilde{A}$-bimodules.
\Ete

Next we prove that $\theta_{A^e}$ is a $(-n)$-degree complex isomorphism in $\mathcal{C}(\tilde{A}, \abimod)$. Analyzing this isomorphism in the derived category, we are led to show that the objects $\RHom_{A^e}(K(A), A^e)$ and $\Hom_{A^e}(K(A), A^e)$ are isomorphic in $\mathcal{D}(\tilde{A},\abimod)$ for any quadratic algebra $A$ over a finite quiver (Proposition \ref{generalphi}). Then the expected result is obtained.

\Bte \label{polynomialstrong}
Let $V$ be a $k$-vector space of finite dimension $n$ and $A=S(V)$ be the symmetric algebra of $V$. Then $A$ is strong $n$-Kc-Calabi-Yau.
\Ete

For the preprojective algebras, we have applied the same strategy in~\cite{bt:kcpa}: state a duality isomorphism $\theta_M$~\cite[Theorem 4.4]{bt:kcpa}, next specialize to $M=A^e$~\cite[Proposition 4.7]{bt:kcpa}, and finally use $\RHom_{A^e}(K(A), A^e) \cong \Hom_{A^e}(K(A), A^e)$ in $\mathcal{D}(\tilde{A},\abimod)$. It would be interesting to implement this strategy for other quadratic algebras, commutative or non commutative.

\setcounter{equation}{0}

\section{A general result on categories of dg modules} \label{generalresult}

\subsection{The derived category associated to a dg algebra and an abelian category} \label{generaldercat}

We follow Yekutieli's book~\cite{yeku:dercat}. In this section, $\kk$ is a nonzero commutative base ring, $B$ is a dg central $\kk$-ring and $\mathcal{M}$ is an abelian category. We shall apply in Subsection \ref{mixedsitu} this general setting to $\mathcal{M}=\almod$ the category of left modules over a central $\kk$-ring $A$, hence the choice of the notation $B$ instead of $A$ used by Yekutieli for the dg algebra. After the Section \ref{generalresult}, $A$ will be a quadratic algebra on a quiver and $B$ will be the Koszul cochain dg algebra of $A$.

Following~\cite{yeku:dercat}, $\mathcal{C}(B,\mathcal{M})$ denotes the dg category of dg left $B$-modules in $\mathcal{M}$ and $\mathcal{C}_{str}(B,\mathcal{M})$ is the associated strict subcategory. The latter is an abelian category. We explicitly describe in Subsection \ref{mixedsitu} these categories when $\mathcal{M}=\almod$. Next $\mathcal{K}(B,\mathcal{M})$ denotes the corresponding homotopy category and the derived category $\mathcal{D}(B,\mathcal{M})$ is obtained from $\mathcal{K}(B,\mathcal{M})$ by inverting quasi-isomorphisms.  One has the canonical functors
$$P: \mathcal{C}_{str}(B,\mathcal{M}) \longrightarrow \mathcal{K}(B,\mathcal{M}),$$
$$Q: \mathcal{K}(B,\mathcal{M}) \longrightarrow \mathcal{D}(B,\mathcal{M}).$$
The localization functor $Q$ is triangulated, meaning that $Q$ respects the triangulated structures of $\mathcal{K}(B,\mathcal{M})$ and $\mathcal{D}(B,\mathcal{M})$.

\subsection{$K$-injective resolutions} \label{existenceright}

Let us recall some definitions and facts on existence of right derived functors~\cite[Section 10.1]{yeku:dercat}. We state them in the general setting considered by Yekutieli. Fix a dg category $\mathcal{C}(B,\mathcal{M})$ as above. An object $N$ in $\mathcal{C}(B,\mathcal{M})$ is said to be \emph{acyclic} if $H^p(N)=0$ for all $p$. An object $I$ in $\mathcal{C}(B,\mathcal{M})$ is said to be \emph{K-injective} if for any acyclic $N$, the object 
$\Hom_{\mathcal{C}(B,\mathcal{M})}(N,I)$ in $\mathcal{C}(\kkmod)$ is acyclic, which is equivalent to say that $\Hom_{\mathcal{K}(B,\mathcal{M})}(N,I)=0$ for every acyclic $N \in \mathcal{K}(B,\mathcal{M})$.

In the $K$-injective terminology as in the $K$-projective terminology, the letter $K$ is referred to the notation used for homotopy categories. $K$-injective (and $K$-projective) resolutions were applied by N. Spaltenstein~\cite{spalt:unbounded} to unbounded complexes in abelian categories and to complexes in categories of sheaves. See also Gelfand-Manin's book~\cite{gm:methods}. The same concepts were emphasized by B. Keller~\cite{keller:derdgcat} in the context of dg modules over dg rings. See~\cite[Remark 10.1.29]{yeku:dercat} for other references.

A \emph{K-injective resolution} of $C\in \mathcal{C}(B,\mathcal{M})$ is a quasi-isomorphism $\rho: C \rightarrow I$ in $\mathcal{C}_{str}(B,\mathcal{M})$, where $I \in \mathcal{C}(B,\mathcal{M})$ is $K$-injective. Quasi-isomorphism means that the morphisms $H^p(\rho)$ in $\mathcal{M}$ are isomorphisms for all $p$. The morphism $P(\rho):C \rightarrow I$ in $\mathcal{K}(B,\mathcal{M})$ is said to be a $K$-injective resolution as well.

As in~\cite[Section 10.1]{yeku:dercat}, $\mathcal{K}(B,\mathcal{M})_{inj}$ denotes the full triangulated subcategory of $\mathcal{K}(B,\mathcal{M})$ on the $K$-injective objects. By definition, $\mathcal{C}(B,\mathcal{M})$ has \emph{enough K-injectives} if any $C \in \mathcal{C}(B,\mathcal{M})$ admits a $K$-injective resolution. We also say that $\mathcal{K}(B,\mathcal{M})$ has enough $K$-injectives. Then, in this case, the restricted localization functor
$$Q:  \mathcal{K}(B,\mathcal{M})_{inj} \rightarrow \mathcal{D}(B,\mathcal{M})$$
will be an equivalence of categories~\cite[Corollary 10.1.15]{yeku:dercat}. 

In this general setting, the existence theorem for right derived functors is the following~\cite[Theorem 10.1.20]{yeku:dercat}.

\Bte \label{rderfunctexist}
Assume that $B$ is a dg central $\kk$-ring and $\mathcal{M}$ is an abelian category such that $\mathcal{K}(B,\mathcal{M})$ has enough $K$-injectives. Let $E$ be a $\kk$-linear triangulated category, and let $F:\mathcal{K}(B,\mathcal{M})\rightarrow E$ be a $\kk$-linear triangulated functor. Then $F$ has a triangulated right derived functor $(RF, \eta^R):\mathcal{D}(B,\mathcal{M})\rightarrow E$ where $\eta^R:F \Rightarrow RF \circ Q$ is an isomorphism of triangulated functors from $\mathcal{K}(B,\mathcal{M})$ to $E$. Furthermore, for every $I \in \mathcal{K}(B,\mathcal{M})_{inj}$, the morphism $\eta^R_I :F(I) \rightarrow RF\circ Q(I)$ in $E$ is an isomorphism.
\Ete

\subsection{$K$-projective resolutions} \label{existenceleft}

Using~\cite[Section 10.2]{yeku:dercat}, there are analogous definitions and facts for $K$-projective objects, $K$-projective resolutions and so on. We just state here the existence theorem for left derived functors~\cite[Theorem 10.2.15]{yeku:dercat}.

\Bte \label{lderfunctexist}
Assume that $B$ is a dg central $\kk$-ring and $\mathcal{M}$ is an abelian category such that $\mathcal{K}(B,\mathcal{M})$ has enough $K$-projectives. Let $E$ be a $\kk$-linear triangulated category, and let $F:\mathcal{K}(B,\mathcal{M})\rightarrow E$ be a $\kk$-linear triangulated functor. Then $F$ has a triangulated left derived functor $(LF, \eta^L):\mathcal{D}(B,\mathcal{M})\rightarrow E$ where $\eta^L:LF \circ Q \Rightarrow F$ is an isomorphism of triangulated functors from $\mathcal{K}(B,\mathcal{M})$ to $E$. Furthermore, for every $P \in \mathcal{K}(B,\mathcal{M})_{prj}$, the morphism $\eta^L_P : LF \circ Q(P) \rightarrow F(P)$ in $E$ is an isomorphism.
\Ete

We shall need another theorem useful for contravariant functors~\cite[Theorem 10.4.8]{yeku:dercat}.

\Bte \label{contrarderfunctexist}
Assume that $B$ is a dg central $\kk$-ring and $\mathcal{M}$ is an abelian category such that $\mathcal{K}(B,\mathcal{M})$ has enough $K$-projectives. Let $E$ be a $\kk$-linear triangulated category, and let $F:\mathcal{K}(B,\mathcal{M})^{op} \rightarrow E$ be a $\kk$-linear triangulated functor. Then $F$ has a triangulated right derived functor $(RF, \eta^R):\mathcal{D}(B,\mathcal{M})^{op} \rightarrow E$ where $\eta^R:F \Rightarrow RF \circ Q$ is an isomorphism of triangulated functors from $\mathcal{K}(B,\mathcal{M})^{op}$ to $E$. Furthermore, for every $P \in \mathcal{K}(B,\mathcal{M})_{prj}$, the morphism $\eta^R_P : F(P)  \rightarrow RF\circ Q(P)$ in $E$ is an isomorphism.
\Ete

Actually, Theorem \ref{rderfunctexist}, Theorem \ref{lderfunctexist} and Theorem \ref{contrarderfunctexist} are stated in~\cite{yeku:dercat} in a bit more general form but the above statements will be sufficient for our applications.

\subsection{Specializations} \label{specializations}

There are two important specializations. If $B=\kk$ is the trivial dg algebra, then $\mathcal{C}(B,\mathcal{M})$ is regarded as the dg category $\mathcal{C}(\mathcal{M})$ of complexes in $\mathcal{M}$. If $\mathcal{M}=\kkmod$, then $\mathcal{C}(B,\mathcal{M})$ is regarded as the dg category $\mathcal{C}(B)$ of dg left $B$-modules. For these specializations, we have the following fundamental results. The references are~\cite[Corollary 11.5.26]{yeku:dercat},~\cite[Corollary 11.3.19]{yeku:dercat},~\cite[Corollary 11.6.30]{yeku:dercat} and~\cite[Corollary 11.4.26]{yeku:dercat} respectively.

\Bte \label{Menoughinj}
If the abelian category $\mathcal{M}$ has enough injectives, the subcategory $\mathcal{C}^{+}(\mathcal{M})$ of bounded below complexes has enough $K$-injectives.
\Ete
\Bte \label{Menoughproj}
If the abelian category $\mathcal{M}$ has enough projectives, the subcategory $\mathcal{C}^{-}(\mathcal{M})$ of bounded above complexes has enough $K$-projectives.
\Ete
\Bte \label{Benoughinj}
If $B$ is any dg central $\kk$-ring, the category $\mathcal{C}(B)$  has enough $K$-injectives.
\Ete
\Bte \label{Benoughproj}
If $B$ is any dg central $\kk$-ring, the category $\mathcal{C}(B)$  has enough $K$-projectives.
\Ete

It is interesting to note that the two specializations intersect. Assume that the dg algebra $B$ is concentrated in degree zero so that the differential is zero. In this case, the dg category $\mathcal{C}(B)$ can be regarded as the dg category $\mathcal{C}(\blmod)$ of complexes of left $B$-modules where $B$ is just considered as a central $\kk$-ring. Then Theorem \ref{Menoughinj} and Theorem \ref{Menoughproj} apply since the abelian category $\blmod$ has enough injectives and enough projectives. However, Theorem \ref{Benoughinj} and Theorem \ref{Benoughproj} say more: the category $\mathcal{C}(\blmod)$ has enough $K$-injectives and enough $K$-projectives \emph{with no boundedness conditions}. In that sense, the two latter theorems are stronger than the two former ones.

The fact that the category $\mathcal{C}(\blmod)$ has enough $K$-injectives and enough $K$-projectives with no boundedness conditions was obtained by N. Spaltenstein~\cite[Theorem C]{spalt:unbounded} and was a motivation for him to introduce $K$-injective and $K$-projective resolutions.

For any dg central $\kk$-ring $B$, the two latter theorems lie on two main theorems of the theory, namely~\cite[Theorem 11.6.24]{yeku:dercat} and~\cite[Theorem 11.4.17]{yeku:dercat} whose proofs are due to B. Keller and are detailed in~\cite{yeku:dercat}. These main theorems are existence theorems of semi-cofree and semi-free resolutions, respectively. Semi-cofree (semi-free) objects in $\mathcal{C}(B)$ are special $K$-injective ($K$-projective) objects that are defined in terms of cofiltrations (filtrations).

\subsection{A mixed specialization $\mathcal{C}(B,\almod)$} \label{mixedsitu}

Our Theorem \ref{mixedsitutheo} below is a new general result in the Yekutieli formalism. This theorem states that $\mathcal{C}(B,\mathcal{M})$ has enough $K$-injectives and enough $K$-projectives in a mixed specialization, that is, besides the specializations $B=\kk$ and $\mathcal{M}=\kkmod$. Let us begin to describe explicitly $\mathcal{C}(B,\mathcal{M})$ in this mixed specialization. The description is drawn from the general definition of $\mathcal{C}(B,\mathcal{M})$ given in~\cite[Chapter 3]{yeku:dercat}.

By default, modules are left modules, graded modules are graded left modules and dg modules are dg left modules. As in~\cite[Chapter 3]{yeku:dercat}, by default, any graded module $M$ is denoted cohomologically $M=\bigoplus_{i\in\zz} M^i$ and the differential of a dg module has degree $+1$. 

In the remainder of the present Section \ref{generalresult}, $A$ is a central $\kk$-ring, $\mathcal{M}=\almod$ is the abelian category of $A$-modules and $B$ is a dg central $\kk$-ring. In the sequel, $\alpha$ and $\beta$ denote elements in $A$ while $f$ and $g$ denote elements in $B$. The product of $f$ and $g$ in $B$ is denoted by $f\times g$ and the unit of $B$ is $1_B$. The differential of $B$ is denoted by $b$.

Let $M=(M^q)_{q \in \zz}$ be a graded $B$-module. This means that $M=(M^q)_{q \in \zz}$ is a graded $\kk$-module and that any $f \in B^p$ acts on $x\in M^q$ by defining an element $f.x$ in $M^{q+p}$, $\kk$-linearly in $f$ and in $x$, such that 
$$f.(g.x))= (f \times g).x, \ \ 1_{B}.x=x.$$
Then $M$ is a graded $B$-module in $\almod$ if moreover each $M^q$ is a module over the $\kk$-ring $A$ such that the actions of $B$ and $A$ on $M$ commute, meaning that the equality
$$\alpha (f.x) = f.(\alpha x)$$
holds for any $\alpha$ in $A$. 

Denote by $\mathcal{G}(B,\almod)$ the category of graded $B$-modules in $\almod$. A morphism in this category is a finite sum of homogeneous graded morphisms. Such an $r$-homogeneous graded morphism $u: M \rightarrow M'$ is $r$-homogeneous graded linear for $B$-actions and each $u: M^q \rightarrow M'^{q+r}$ is a morphism of $A$-modules. The subcategory $\mathcal{G}_{str}(B,\almod)$ whose morphisms are strict, i.e., of degree 0, is abelian; $\mathcal{G}_{str}(B,\almod)$ is called the strict category of graded $B$-modules in $\almod$. 

Up to now, we have just used the graded central $\kk$-ring $B$. Taking into account the differential $b$ of $B$, that is, the whole structure of dg central $\kk$-ring of $B$, we are ready to define the dg category $\mathcal{C}(B,\almod)$ of dg $B$-modules in $\almod$. An object $M=(M^q)_{q\in \zz}$ in $\mathcal{C}(B,\almod)$ is a graded $B$-module in $\almod$ endowed with a differential $d_M: M^q \rightarrow M^{q+1}$ of $A$-modules, satisfying 
\begin{equation} \label{axiom}
  d_M(f.x)=b(f).x + (-1)^p f.d_M(x)
 \end{equation}
for any $f\in B^p$ and $x$ in $M$. The forgetful functor $\Und: \mathcal{C}(B,\almod) \rightarrow \mathcal{G}(B,\almod)$ that forgets the differentials of $B$ and of the dg $B$-modules is a fully faithful graded functor~\cite[Proposition 3.8.18]{yeku:dercat}.

By definition, for two objects $M$ and $M'$ in $\mathcal{C}(B,\almod)$, one has
$$\Hom_{\mathcal{C}(B,\almod)}(M,M')= \Hom_{\mathcal{G}(B,\almod)}(M,M'),$$
and this graded space has a differential $d$ defined by
\begin{equation} \label{diffofu}
d(u)=d_{M'} \circ u - (-1)^r u \circ d_M 
\end{equation}
where $u : M \rightarrow M'$ is $r$-homogeneous. Then $(\Hom_{\mathcal{C}(B,\almod)}(M,M'),d)$ is a complex in $\kkmod$, i.e., a dg $\kk$-module.

The subcategory $\mathcal{C}_{str}(B,\almod)$ of $\mathcal{C}(B,\almod)$ has the same objects, with morphism spaces 
$$\Hom_{\mathcal{C}_{str}(B,\almod)}(M,M')= Z^0 (\Hom_{\mathcal{C}(B,\almod)}(M,M')),$$
that is, the morphisms from $M$ to $M'$ are of degree 0 \emph{and} are morphisms of complexes. Warning: these morphisms are still called strict morphisms. The abelian category $\mathcal{C}_{str}(B,\almod)$ is called the strict category of dg $B$-modules in $\almod$. The forgetful functor $\Und: \mathcal{C}_{str}(B,\almod) \rightarrow \mathcal{G}_{str}(B,\almod)$ is faithfully exact.

The homotopy category $\mathcal{K}(B,\almod)$ has the same objects as $\mathcal{C}(B,\almod)$, with morphism spaces
$$\Hom_{\mathcal{K}(B,\almod)}(M,M')= H^0 (\Hom_{\mathcal{C}(B,\almod)}(M,M')),$$
that is, the morphisms from $M$ to $M'$ are the homotopy classes of strict morphisms. It is a triangulated category. The canonical functor
$$P: \mathcal{C}_{str}(B,\almod) \longrightarrow \mathcal{K}(B,\almod)$$
sends strict morphisms to their homotopy classes.

From $\mathcal{K}(B,\almod)$, the derived category $\mathcal{D}(B,\almod)$ is constructed by inverting quasi-isomorphisms. This is a triangulated category and the localization functor
$$Q: \mathcal{K}(B,\almod) \longrightarrow \mathcal{D}(B,\almod)$$
is triangulated.

\subsection{A dg functor $\Phi$} \label{functorPhi}

Our Theorem \ref{mixedsitutheo} is based on a functor $\Phi$ defined now. In all the sequel, by an algebra, we mean a central $\kk$-ring. The unadorned tensor symbol $\otimes$ means $\otimes_{\kk}$.

The algebra $A$ is endowed with a structure of trivial dg algebra by imposing the grading $A=A^0$ and the zero differential. Therefore, $A \otimes B$ is a dg algebra~\cite[Example 3.3.10]{yeku:dercat}. The product is given by
$$(\alpha \otimes f)\times (\beta \otimes g)=(\alpha \beta)\otimes (f\times g)$$
for $f, g$ in $B$ and $\alpha, \beta$ in $A$. The unit is $1_A \otimes 1_{B}$. The differential $\check{b}$ is given by
$$\check{b}(\alpha \otimes f)=\alpha \otimes b(f).$$

Our aim is to define a dg functor
\begin{equation} \label{deffunctorF}
\Phi : \mathcal{C}(B, \almod) \longrightarrow \mathcal{C}(A \otimes B).
\end{equation} 
Let $M=(M^q)_{q \in \zz}$ be a dg $B$-module in $\almod$. We define actions of $A \otimes B$ on $x \in M$ by 
\begin{equation}
(\alpha \otimes f).x  = \alpha(f.x)  = f.(\alpha x).
\end{equation}
The equality on the right is the commutation of $B$ and $A$ acting on $M$. Clearly
\begin{align}
(\alpha \otimes f).((\beta \otimes g).x)  & = ((\alpha \otimes f)\times(\beta \otimes g)).x \, ,\\
(1_A\otimes 1_{B}).x & = x  \, .  
\end{align}
So $M$ is a graded $(A \otimes B)$-module.

Let $f\in B^p$. On one hand, one has
$$d_M((\alpha \otimes f).x)= d_M(f.(\alpha x))=b(f).(\alpha x))+(-1)^p f.d_M( \alpha x)= (\alpha \otimes b(f)).x + (-1)^p f.(\alpha d_M(x))$$
since $d_M$ is $A$-linear. On the other hand,
$$\check{b} (\alpha \otimes f).x + (-1)^p (\alpha \otimes f).d_M(x)=(\alpha \otimes b(f)).x + (-1)^p f.(\alpha d_M(x))$$
and we obtain 
\begin{equation} \label{diffacting1}
d_M((\alpha \otimes f).x) = \check{b}(\alpha \otimes f).x + (-1)^p (\alpha \otimes f).d_M(x).  
\end{equation}
Thus $M$ is a dg $(A \otimes B)$-module, defining an object $\Phi(M)$ in $\mathcal{C}(A \otimes B)$. Note that the differential of the object $\Phi(M)$ is just $d_M$ considered as a $\kk$-linear map.

Define now $\Phi$ on morphisms. Let $u: M \rightarrow M'$ be an $r$-homogeneous graded morphism for $B$-actions such that each $u:M^q \rightarrow M'^{q+r}$ is a morphism of $A$-modules. The equality
$$ u((\alpha \otimes f).x)=u(f.(\alpha x))=(-1)^{pr}f.u(\alpha x)=(-1)^{pr}f.(\alpha u(x))=(-1)^{pr}(\alpha \otimes f).u(x)$$
for $f\in B^p$ shows that $u: M \rightarrow M'$ becomes a morphism in $\mathcal{C}(A \otimes B)$ from $\Phi(M)$ to $\Phi(M')$, denoted by $\Phi(u)$. The functoriality properties $\Phi(\id_C)=\id_{\Phi(C)} $ and $\Phi(v \circ u)=\Phi(v)\circ \Phi(u)$ hold.
Consequently $\Phi$ is a functor. Note that we can extend $\Phi$ to any $A$-linear map from $M$ to $M'$ in an obvious manner, so that we can write $\Phi(d_M)=d_{\Phi(M)}$.

Let us prove that the functor $\Phi$ is a dg functor~\cite[Definition 3.5.1]{yeku:dercat}. Recall from Equality (\ref{diffofu}) that
$$d(u)=d_{M'} \circ u - (-1)^r u \circ d_M\, .$$
By functoriality $\Phi(d(u))=\Phi(d_{M'}) \circ \Phi(u) - (-1)^r \Phi(u) \circ \Phi(d_M)$, implying that $\Phi(d(u))=d(\Phi(u))$, as we want. In conclusion, we obtain the following.

\Bpo \label{Ffuncto}
As defined above, $\Phi: \mathcal{C}(B, \almod) \longrightarrow \mathcal{C}(A \otimes B)$ is a dg functor.
\Epo

\subsection{A dg functor $\Psi$, inverse to $\Phi$} \label{functorPsi}

Let $(M,d_M)$ be a dg $(A \otimes B)$-module. For $f$ in $B^p$ and $\alpha\in A$, we define the following actions
on $x \in M$ by
\begin{align}
f.x  & = (1_A \otimes f).x \\
\alpha x  & = (\alpha \otimes 1_{B}).x .
\end{align}
For these actions, $M$ becomes a graded $B$-module in $\almod$. Moreover, one has
$$d_M(f.x)=d_M((1_A \otimes f).x)=\check{b}(1_A\otimes f).x + (-1)^p(1_A \otimes f).d_M(x),$$
and using $\check{b}(1_A\otimes f)=1_A\otimes b(f)$, we obtain
$$d_M(f.x)=b(f).x + (-1)^p f.d_M(x).$$
It is immediate that $d_M$ is $A$-linear. Therefore, $M$ becomes a dg $B$-module in $\almod$, denoted by $\Psi(M)$.

Let $M$ and $M'$ be dg $(A \otimes B)$-modules. Let $u: M \rightarrow M'$ be an $r$-homogeneous graded morphism for $(A \otimes B)$-actions. Then it is direct to prove that $u$ is a graded $B$-module morphism and an $A$-module morphism, defining $\Psi(u):\Psi(M)\rightarrow \Psi(M')$. Clearly, we get a functor
\begin{equation} \label{deffunctorG}
\Psi: \mathcal{C}(A \otimes B) \rightarrow \mathcal{C}(B, \almod).
\end{equation}

Here again, we can extend $\Psi$ to any $A$-linear map from $M$ to $M'$, so that we can write $\Psi(d_M)=d_{\Psi(M)}$ and deduce $\Psi(d(u))=d(\Psi(u))$. Thus $\Psi$ is a dg functor. Since $\Psi\circ \Phi$ and $\Phi \circ \Psi$ are identity functors, we have proved the following.

\Bpo \label{fundatheo2}
The functor $\Phi: \mathcal{C}(B, \almod) \longrightarrow \mathcal{C}(A \otimes B)$ is an isomorphism of dg categories.
\Epo

\subsection{A result in the mixed specialization} \label{resultmixed}

Since $\Phi(M)$ and $M$ have the same $\kk$-linear differential, one has the commutative diagram
\begin{eqnarray} \label{HPhi}
 \mathcal{C}_{str}(B,\almod) \ \ \ \stackrel{\Phi}{\longrightarrow} & \mathcal{C}_{str}(A\otimes B) \nonumber  \\
\downarrow H  \ \ \ \ \ \ \ \ \ \ \   &  \downarrow H \\
\mathcal{G}_{str}(\almod) \ \ \ \stackrel{\Und}{\longrightarrow} & \mathcal{G}_{str}(\kkmod) \nonumber
\end{eqnarray}
where $\Und$ is a forgetful functor.
\Bpo \label{stabilityPhi}
(i) The object $M$ is acyclic in $\mathcal{C}(B,\almod)$ if and only if $\Phi(M)$ is acyclic in $\mathcal{C}(A\otimes B)$.

(ii) A strict morphism $u$ in $\mathcal{C}(B, \almod)$ is a quasi-isomorphism if and only if $\Phi(u)$ is a quasi-isomorphism in $\mathcal{C}(A \otimes B)$.

(iii) An object $M$ is $K$-injective ($K$-projective) in $\mathcal{C}(B,\almod)$ if and only if $\Phi(M)$ is $K$-injective ($K$-projective) in $\mathcal{C}(A\otimes B)$.
\Epo
\Bdm
(i) and (ii) are clear from the commutative diagram ($\ref{HPhi}$). To deduce (iii) from (i), remark that the dg isomorphism $\Phi$ induces an isomorphism
$$\Hom_{\mathcal{C}(B,\almod)}(M,M') \longrightarrow \Hom_{\mathcal{C}(A\otimes B)}(\Phi(M),\Phi(M'))$$
in $\mathcal{C}_{str}(\kkmod)$, so that the complex $\Hom_{\mathcal{C}(B,\almod)}(M,M')$ is acyclic if and only if the complex $\Hom_{\mathcal{C}(A\otimes B)}(\Phi(M),\Phi(M'))$ is acyclic. Therefore, an object $M$ is $K$-injective in $\mathcal{C}(B,\almod)$ if and only if the complex $\Hom_{\mathcal{C}(A\otimes B)}(\Phi(N),\Phi(M))$ is acyclic for any acyclic object $N$ in $\mathcal{C}(B,\almod)$. Since $\Phi$ is an isomorphism, any acyclic object $N'$ in $\mathcal{C}(A\otimes B)$ can be written $N'=\Phi(N)$, and $N$ is acyclic in $\mathcal{C}(B,\almod)$ by (i). Thus $M$ is $K$-injective in $\mathcal{C}(B,\almod)$ if and only if $\Phi(M)$ is $K$-injective in $\mathcal{C}(A\otimes B)$. There is an analogous proof in the $K$-projective case.
\Edm
\\

As a consequence of (ii) and (iii), $\rho: M \rightarrow I$ is a $K$-injective resolution in $\mathcal{C}_{str}(A \otimes B)$ if and only if $\Phi^{-1}(\rho): \Phi^{-1}(M) \rightarrow \Phi^{-1}(I)$ is a $K$-injective resolution in $\mathcal{C}_{str}(B, \almod)$. There is an analogous statement in the $K$-projective case. Therefore, Theorem \ref{Benoughinj} and Theorem \ref{Benoughproj} imply the main result of this section.

\Bte \label{mixedsitutheo}
Assume that $A$ is a central $\kk$-ring and $B$ is a dg central $\kk$-ring. The category $\mathcal{C}(B, \almod)$ has enough $K$-injectives and enough $K$-projectives.
\Ete
\Bdm
Let $C$ be an object in $\mathcal{C}_{str}(B, \almod)$. Then $M=\Phi(C)$ has a $K$-injective resolution $\rho: M \rightarrow I$ in $\mathcal{C}_{str}(A \otimes B)$ by Theorem \ref{Benoughinj}, so that $\Phi^{-1}(\rho): C \rightarrow \Phi^{-1}(I)$ is a $K$-injective resolution in $\mathcal{C}_{str}(B, \almod)$. Use Theorem \ref{Benoughproj} for an analogous proof in the $K$-projective case.
\Edm
\\

This theorem is a generalization of Spaltenstein's theorem claiming that $\mathcal{C}(\almod)$ has enough $K$-injectives and enough $K$-projectives (see comments after Theorem \ref{Benoughproj}).

\subsection{A bimodule version of $\mathcal{C}(B,\almod)$} \label{bimversion}

Let us keep the assumptions of the mixed specialization as presented in Subsection \ref{mixedsitu}, that is, $A$ is a central $\kk$-ring and $B$ is a dg central $\kk$-ring. Let us introduce the opposite central $\kk$-ring $A^{op}$ and the opposite dg central $\kk$-ring $B^{op}$~\cite[Definition 3.3.13]{yeku:dercat}. The multiplication $\times^{op}$ of $B^{op}$ is reversed and twisted by Koszul signs: $f\times^{op} g:= (-1)^{pq} g\times f$ for $f\in B^p$ and $g\in B^q$, while the multiplication of $A^{op}$ is just reversed. The differential of $B^{op}$ is equal to the differential of $B$.

Next, we define the central $\kk$-ring $A^e:= A \otimes A^{op}$ and the dg central $\kk$-ring $B^e:=B \otimes B^{op}$, called the \emph{enveloping} ring (dg ring) of $A$ (of $B$). The multiplication of $B^e$ is given by
\begin{equation} \label{multenv}
(f_1\otimes g_1)\times (f_2 \otimes g_2):=(-1)^{q_1p_2} (f_1\times f_2) \otimes (g_1 \times^{op} g_2)
\end{equation}
for any $f_1 \in B$, $g_1 \in B^{q_1}$, $f_2 \in B^{p_2}$ and $g_2 \in B$. The differential $\tilde{b}$ of $B^e$ is given by
\begin{equation} \label{diffenv}
\tilde{b}(f\otimes g):=b(f) \otimes g + (-1)^p f \otimes b(g)
\end{equation}
where $f \in B^p$ and $g \in B$~\cite[Example 3.3.10]{yeku:dercat}. Theorem \ref{mixedsitutheo} applies to $A^e$ and $B^e$, which provides an enveloping version of this theorem.

\Bcr \label{envversiontheo}
Assume that $A$ is a central $\kk$-ring and $B$ is a dg central $\kk$-ring. The category $\mathcal{C}(B^e, \aemod)$ has enough $K$-injectives and enough $K$-projectives.
\Ecr

As it is well-known~\cite{ce:homo, weib:homo}, the main feature of $A^e$ is that an $A^e$-module $M$ is the same thing as a $\kk$-central $A$-bimodule via the actions
\begin{equation} \label{bimodact1}
\alpha x \beta:= (\alpha \otimes \beta) x
\end{equation}
for any $\alpha \in A$, $\beta \in A$ and $x \in M$. 

Similarly, a graded $B^e$-module $M=(M^q)_{q \in \zz}$ is the same thing as a $\kk$-central graded $B$-bimodule via the actions 
\begin{equation} \label{bimodact2}
f.x.g:=(-1)^{rq} (f\otimes g).x
\end{equation}
for any $f\in B^p$, $g\in B^q$ and $x\in M^r$.

Therefore, a graded $B^e$-module $M=(M^q)_{q \in \zz}$ in the abelian category $\aemod$ is the same thing as a $\kk$-central graded $B$-bimodule in the abelian category $\abimod$ of $A$-bimodules, meaning that
$$\alpha (f.x.g)\beta = f.(\alpha x \beta).g$$
holds for any $\alpha$ and $\beta$ in $A$.

If moreover $(M,d_M)$ is a dg $B^e$-module in $\aemod$, each $d_M:M^q \rightarrow M^{q+1}$ is a morphism of $A$-bimodules for the actions (\ref{bimodact1}), satisfying 
\begin{equation} \label{diffproperties}
  d_M(f.x)=b(f).x + (-1)^p f.d_M(x), \ \ d_M(x.f)=d_M(x).f + (-1)^q x.b(f)
 \end{equation}
for any $f\in B^p$ and $x$ in $M^q$. We leave the verifications to the reader. Then we say that $M$ is a $\kk$-central dg $B$-bimodule in $\abimod$. Doing so, the dg category $\mathcal{C}(B^e, \aemod)$ is identified to the dg category denoted by $\mathcal{C}_{bim}(B,\abimod)$ of $\kk$-central dg $B$-bimodules in $\abimod$. Consequently, Corollary \ref{envversiontheo} is translated in a bimodule version as follows.

\Bcr \label{bimversiontheo}
Assume that $A$ is a central $\kk$-ring and $B$ is a dg central $\kk$-ring. The category $\mathcal{C}_{bim}(B, \abimod)$ has enough $K$-injectives and enough $K$-projectives.
\Ecr

\setcounter{equation}{0}

\section{The derived Koszul calculus of a quadratic algebra} \label{derkoscalA}

\subsection{The Koszul calculus of a quadratic algebra} \label{koscalA}

Throughout the remainder of the paper, we follow the setup of~\cite{bt:kcpa}. Let $\mathcal{Q}=(\mathcal{Q}_0,\mathcal{Q}_1)$ be a quiver such that the set $\mathcal{Q}_0$ of vertices and the set $\mathcal{Q}_1$ of arrows are finite. We assume that $\mathcal{Q}_0 \neq \emptyset$. Let $\ff$ be a field, playing the role of the base ring $\kk$ of the general setting presented in Section $\ref{generalresult}$. The unadorned tensor symbol $\otimes$ means $\otimes_{\ff}$.

The vertex space $k=\ff \mathcal{Q}_0$ becomes a (nonzero) commutative ring isomorphic to the ring $\ff^{|\mathcal{Q}_0|}$ by associating with $\mathcal{Q}_0$ a complete set of orthogonal idempotents $\{e_i\,;\,i\in \mathcal{Q}_0\}$. For each arrow $\alpha \in \mathcal{Q}_1$, denote its source vertex by $\mo(\alpha)$ and its target vertex by $\mathfrak{t}(\alpha)$. The arrow space $V=\ff \mathcal{Q}_1$ is a $k$-bimodule for the following actions: $e_j\alpha e_i$ is equal to zero if $i\neq \mathfrak{s}(\alpha)$ or $j\neq \mathfrak{t}(\alpha)$, and is equal to $\alpha$ if $i= \mathfrak{s}(\alpha)$ and $j= \mathfrak{t}(\alpha)$. 

The tensor $k$-algebra $T_k(V)$ of the $k$-bimodule $V$ is identified to the path algebra $\mathbb{F}\mathcal{Q}$ of the quiver $\mathcal{Q}$. Fixing a $k$-subbimodule $R$ of $V\otimes_kV\cong \mathbb{F}\mathcal{Q}_2$, the associative algebra
$$A=T_k(V)/(R)\cong \mathbb{F}\mathcal{Q}/(R)$$
where $(R)$ denotes the two-sided ideal generated by $R$, is called a \emph{quadratic $k$-algebra over the finite quiver $\mathcal{Q}$}.

The algebra $A$ is $\mathbb{F}$-central, but is not $k$-central in general. All the bimodules considered in this paper are assumed to be $\ff$-central. Setting $A^e=A\otimes A^{op}$, an $A$-bimodule can be viewed as a left (or right) $A^e$-module, as usual. The algebra $A$ is graded for the path length. The degree $m$ component of $A$ is denoted by $A_m$. Clearly, $A_0\cong k$ and $A_1\cong V$.

The \emph{Koszul complex} $K(A)$ is the subcomplex of the bar resolution $B(A)$ defined by the $A$-subbimodules $A\otimes_k W_p \otimes_k A$ of $A\otimes_k A^{\otimes_k p} \otimes_k A$, where $W_0=k$, $W_1=V$, $W_2=R$ and, for $p\pgq 3$,
\begin{equation} \label{definitionw}
  W_{p}=\bigcap_{i+2+j=p}V^{\otimes_k i}\otimes_k R\otimes_k V^{\otimes_k j}
\end{equation}
with $i\geq 0$ and $j\geq 0$. Here $W_p$ is considered as a $k$-subbimodule of $V^{\otimes_k p} \subseteq A^{\otimes_k p}$. The differential $d$ of the chain complex $K(A)$ is defined on $K(A)_p=A\otimes_k W_p \otimes_k A$ by
\begin{equation} \label{definitiond}
d(\alpha \otimes_k x_1 \ldots x_{p} \otimes_k \beta) =\alpha x_1\otimes_k x_2 \ldots x_{p}\otimes_k \beta +(-1)^p \alpha \otimes_k x_1 \ldots x_{p-1}\otimes_k x_{p} \beta,
\end{equation}
for $\alpha$, $\beta$ in $A$ and $x_1 \ldots x_{p}$ in $W_p$. The $A$-bimodules $A\otimes_k W_p \otimes_k A$ are finitely generated. So
$$K(A): \ \ \ \ \ \cdots \stackrel{d}{\longrightarrow} A \otimes_k R \otimes_k A \stackrel{d}{\longrightarrow} A \otimes_k V \otimes_k A \stackrel{d}{\longrightarrow} A \otimes_k A \longrightarrow 0.$$

In this paper, we follow~\cite{bt:kcpa} for the notation of elements of $W_p$. Let us recall this notation. As in (\ref{definitiond}), an arbitrary element of $W_p$ is denoted by a product $x_1 \ldots x_{p}$ thought of as a sum of such products, where $x_1, \ldots , x_p$ are in $V$. Moreover, regarding $W_p$ as a subspace of $V^{\otimes_k q}\otimes_k W_r \otimes_k V^{\otimes_k s}$ with $q+r+s=p$, the 
element $x_1 \ldots  x_p$ viewed in $V^{\otimes_k q}\otimes_k W_r \otimes_k V^{\otimes _k s}$ will be denoted by the \emph{same} notation, meaning that 
the product $x_{q+1} \ldots x_{q+r}$ represents an element of $W_r$ and the other $x_i$'s are viewed as arbitrary in $V$.

 As the complex $B(A)$, $K(A)$ is a chain complex of projective $A$-bimodules ending in degree 0 with $A \otimes_k A$. The bar resolution $B(A)$ of $A$ is defined by the multiplication $\mu: A \otimes_k A \rightarrow A$ of $A$. In general, $K(A)$ is not a subresolution of $B(A)$. When it is a subresolution, $A$ is said to be \emph{Koszul}. If $A$ is Koszul, it is better to work with $K(A)$ since $K(A)$ is then a minimal projective resolution of $A$~\cite{bt:kcpa}.
\Brm
This definition of being Koszul is equivalent to the standard one expressed in terms of resolutions of $k$ by projective left $A$-modules~\cite{vdb:nch}. See also~\cite{rbnm:kogo} for a proof of this equivalence in the more general context of $N$-homogeneous algebras. 
\Erm

Let us continue the setup. We associate to $A$ the dg central $\ff$-ring
$$\tilde{A}=(\Hom_{A^e}(K(A), A), b_K, \cupk).$$
We simply say that $\tilde{A}$ is a dg algebra. For the definitions of $b_K$ and $\cupk$, see Equalities (\ref{defb_K}) and (\ref{defcupK}) below. The elements of $\tilde{A}$ should be thought of as differentials on $A$. If $A$ is a polynomial algebra in $n$ variables, the algebra $\tilde{A}$ is isomorphic to the exterior algebra $\Omega_A^{\ast}$ constructed on the K\"{a}hler differentials of $A$~\cite{weib:homo}, as we shall see it in Theorem \ref{dualityiso}.

\Bdf \label{defAtilde}
The dg algebra $\tilde{A}$ is called the Koszul cochain dg algebra of $A$.
\Edf

As in~\cite{bt:kcpa}, for any $A$-bimodule $M$, the complex $(\Hom_{A^e}(K(A), M),b_K)$ whose homology defines the Koszul cohomology space $\HK^{\bullet}(A,M)$, is identified to $(\Hom_{k^e}(W_{\bullet},M),b_K)$ where $k^e=k\otimes k$. Similarly, the complex $(M\otimes_{A^e} K(A),b^K)$ whose homology defines the Koszul homology space $\HK_{\bullet}(A,M)$, is identified to $(M\otimes_{k^e} W_{\bullet}, b^K)$. Since $A$ and $M$ are $\ff$-central, the dg algebra $\tilde{A}$ is $\ff$-central, and the spaces $\Hom_{k^e}(W_{\bullet},M)$ and $M\otimes_{k^e} W_{\bullet}$ are $\ff$-central.

\Bdf \label{Koszulspaces}
For $p\geq 0$, the space $Hom_{k^e}(W_p,M)$ ($M\otimes_{k^e} W_p$) is called the space of Koszul $p$-cochains ($p$-chains) of $A$ with coefficients in the $A$-bimodule $M$.
\Edf
\Bpo \label{fundamcomplexes}
For any $A$-bimodule $M$, $(\Hom _{k^e}(W_{\bullet},M),b_K)$ and  $(M\otimes_{k^e} W_{\bullet},b^K)$ are dg bimodules over the dg algebra $\tilde{A}$ for the actions given by the Koszul cup product $\underset{K}{\smile}$ and the Koszul cap product $\underset{K}{\frown}$ respectively. 
\Epo
For a proof and the definitions of the Koszul cup and cap products, see~\cite[Subsection 2.3]{bt:kcpa}. To illustrate these definitions, we give here the definitions of $b_K$ and $\underset{K}{\smile}$. For $f\in \Hom _{k^e}(W_p,M)$, $b_K(f)\in \Hom _{k^e}(W_{p+1},M)$ is defined by 
\begin{equation} \label{defb_K}
  b_K(f)( x_1 \ldots x_{p+1}) =f(x_1\ldots x_{p})\, x_{p+1} -(-1)^p x_1\, f(x_2 \ldots x_{p+1})
  \end{equation}
where $x_1 \ldots x_{p+1} \in W_{p+1}$. For $f\in \Hom _{k^e}(W_p,M)$ and $g\in \Hom _{k^e}(W_q,N)$, we define $f\underset{K}{\smile}g \in \Hom _{k^e}(W_{p+q},M\otimes_A N)$ by
\begin{equation} \label{defcupK}
  (f\underset{K}{\smile} g) (x_1 \ldots x_{p+q}) = (-1)^{pq} f(x_1 \ldots x_p)\otimes_A \, g(x_{p+1} \ldots  x_{p+q})
\end{equation}
where $x_1 \ldots x_{p+q} \in W_{p+q}$.

Let us recall the \emph{fundamental formulas of Koszul calculus}~\cite[Subsection 2.4]{bt:kcpa}. These formulas express the differentials of the complexes $(\Hom _{k^e}(W_{\bullet},M),b_K)$ and  $(M\otimes_{k^e} W_{\bullet},b^K)$ in terms of the Koszul cup and cap products. See ~\cite[Subsection 2.4]{bt:kcpa} for the definitions of the graded commutators $[-, -]_{\underset{K}{\smile}}$ and $[-, -]_{\underset{K}{\frown}}$. Recall that $\me_A: V \rightarrow A$ defined by $\me_A(x)=x$ for all $x\in V$ is a Koszul 1-cocycle, called the \emph{fundamental Koszul 1-cocycle} of $A$.

\Bpo \label{fundamental}
For any Koszul cochain $f$ and any Koszul chain $z$ with coefficients in an $A$-bimodule $M$, we have
\begin{align}
\label{fundacoho}
 b_K(f)&=- [\me_A, f]_{\underset{K}{\smile}},\\
 \label{fundaho}
 b^K(z)&=-[\me_A, z]_{\underset{K}{\frown}}. 
  \end{align}
\Epo

As a consequence of Proposition \ref{fundamcomplexes}, for any $A$-bimodule $M$, $\HK^{\bullet}(A,M)$ and $\HK_{\bullet}(A,M)$ are graded bimodules over the graded $\ff$-algebra $\HK^{\bullet}(A)$, forming the \emph{Koszul calculus} of $A$ in the terminology of~\cite{bt:kcpa}. The Koszul calculus of $A$ only depends on the associative algebra structure of $A$~\cite[Subsection 2.7]{bt:kcpa}. Since $K(A)$ is a complex of projective $A$-bimodules, $M\mapsto \HK^{\bullet}(A,M)$ and $M\mapsto \HK_{\bullet}(A,M)$ define $\delta$-functors from $\abimod$ to $\fmod$, that is, any short exact sequence of $A$-bimodules naturally gives rise to a long exact sequence
in Koszul cohomology and in Koszul homology~\cite[Chapter 2]{weib:homo}.

When $A$ is Koszul, the graded algebra $\HK^{\bullet}(A)$ is isomorphic to the Hochschild cohomology algebra $\HH^{\bullet}(A)$. Moreover $\HK^{\bullet}(A,M)$ and $\HK_{\bullet}(A,M)$ are isomorphic as $\HH^{\bullet}(A)$-bimodules to the Hochschild cohomology space $\HH^{\bullet}(A,M)$ and to the Hochschild homology space $\HH_{\bullet}(A,M)$ respectively, so that the Koszul calculus coincides with the Hochschild calculus when $A$ is Koszul.

A family of non Koszul algebras is given by the preprojective algebras of ADE Dynkin types excluding types A$_1$ and A$_2$. For this family of examples, Koszul calculus and Hochschild calculus are different w.r.t. the Calabi-Yau property~\cite{bt:kcpa}. See below Subsection \ref{definingstrongKcCY}. Moreover, the explicit computations obtained in~\cite{bt:kcpa} for this family of examples show that the Koszul calculus provides more cohomological invariants than the Hochschild calculus, except in type $\mathrm{E}_8$ with $\car(\ff)=2$ where the invariants are the same. 
\\

Let us end this setup by some elementary calculations. The space $\tilde{A}^p$ of Koszul $p$-cochains with coefficients in $A$ is the space of the $k$-bimodule morphisms from $W_p$ to $A$. In particular, the space $\tilde{A}^0$ of $0$-cochains is isomorphic to $\bigoplus_{i \in \mathcal{Q}_0} e_iAe_i$ which is a $k$-subalgebra of $A$ containing the unit $1_A=\sum_{i \in \mathcal{Q}_0} e_i$ of $A$. So the unit $1_{\tilde{A}}$ of $\tilde{A}$ is identified to the unit of $A$. Note that if $f:W_p \rightarrow M$ is a Koszul $p$-cochain with coefficients in an $A$-bimodule $M$ and if $a=\sum_{i \in \mathcal{Q}_0} e_ia_ie_i$ is the $0$-cochain $k \rightarrow A, 1 \mapsto a$, then $a\cupk f$ is the $p$-cochain $x_1 \ldots x_{p} \mapsto a\,f(x_1 \ldots x_{p})$, equivalently for $f \cupk a$. Similarly, if $z=m\otimes_{k^e}x_1 \ldots x_{q}$ is a Koszul $q$-chain with coefficients in $M$, then $a\capk z=(a\,m)\otimes_{k^e}x_1 \ldots x_{q}$, equivalently for $z \capk a$.

\Brm \label{hyper}
As graded $\HK^{\bullet}(A)$-bimodules, $\HK^{\bullet}(A,M)$ and $\HK_{\bullet}(A,M)$ are isomorphic to the Hochschild hypercohomology $\mathbb{H}\mathbb{H}^{\bullet}(A, \Hom _A (K(A),M))$ and to the Hochschild hyperhomology $\mathbb{H}\mathbb{H}_{\bullet}(A, M \otimes_A K(A))$ respectively~\cite[Section 2]{bls:kocal}~\cite[Section 9]{bs:cupcap}. When $A$ is Koszul, $K(A)\cong A$ in $\mathcal{D}(\abimod)$ and we recover that $\HK^{\bullet}(A,M)$ and $\HK_{\bullet}(A,M)$ are isomorphic as $\HH^{\bullet}(A)$-bimodules to $\HH^{\bullet}(A,M)$ and to $\HH_{\bullet}(A,M)$ respectively.
\Erm

\subsection{A dg category of dg $\tilde{A}$-bimodules} \label{dgbimod}

Let us fix a quadratic $k$-algebra $A$ over $\mathcal{Q}$ and consider the Koszul cochain dg algebra $\tilde{A}$. Subsection \ref{bimversion} allows us to introduce the dg category $\mathcal{C}_{bim}(\tilde{A}, \abimod)$.

For simplicity, $\mathcal{C}_{bim}(\tilde{A}, \abimod)$ is denoted by $\mathcal{C}(\tilde{A}, \abimod)$. In other words, \emph{throughout the remainder of the paper}, $\mathcal{C}(\tilde{A}, \abimod)$ denotes the dg category of dg $\tilde{A}$-bimodules in $\abimod$. Corollary \ref{bimversiontheo} is exemplified in the following fundamental result.

\Bpo \label{fundaresult}
As defined above, the category $\mathcal{C}(\tilde{A}, \abimod)$ has enough $K$-injectives and enough $K$-projectives.
\Epo

Let us present in some details the dg category $\mathcal{C}(\tilde{A}, \abimod)$, the associated strict category $\mathcal{C}_{str}(\tilde{A}, \abimod)$, the corresponding homotopy category $\mathcal{K}(\tilde{A},\abimod)$ and derived category $\mathcal{D}(\tilde{A},\abimod)$. We follow the Yekutieli formalism developed in Subsection \ref{mixedsitu}, applied now in our bimodule situation. 

Let $C=(C^q)_{q \in \zz}$ be a graded $\tilde{A}$-bimodule. This means that any Koszul $p$-cochain $f$ with coefficients in $A$ acts on $x\in C^q$ by defining elements $f.x$ and $x.f$ in $C^{q+p}$, $\ff$-linearly in $f$ and in $x$, such that 
$$f.(g.x))= (f \cupk g).x, \ (x.f).g=x.(f \cupk g), \ (f.x).g=f.(x.g), \ 1_{\tilde{A}}.x=x.1_{\tilde{A}}=x.$$
Then $C$ is a graded $\tilde{A}$-bimodule in the abelian category $\abimod$ if moreover, each $C^q$ is an $A$-bimodule such that the actions of $\tilde{A}$ and $A$ on $C$ commute, meaning that the equality 
$$\alpha (f.x.g)\beta = f.(\alpha x \beta).g$$
holds for any $\alpha$ and $\beta$ in $A$. Note that if $f:k \rightarrow A, 1 \mapsto a$ is a $0$-cochain where $a=\sum_{i \in \mathcal{Q}_0} e_ia_ie_i$, then $f.x \neq ax$ in general, for example when $C=K(A)$ as we shall see in Subsection  \ref{K(A)dgbimod}.

Denote by $\mathcal{G}(\tilde{A},\abimod)$ the category of graded $\tilde{A}$-bimodules in $\abimod$. A morphism in this category is a finite sum of homogeneous graded morphisms. Such an $r$-homogeneous graded morphism $u: C \rightarrow C'$ is $r$-homogeneous graded for $\tilde{A}$-actions and each $u: C^q \rightarrow C'^{q+r}$ is a morphism of $A$-bimodules. The subcategory $\mathcal{G}_{str}(\tilde{A},\abimod)$ whose morphisms are strict, i.e., of degree 0, is abelian.

Taking into account the differential $b_K$ of $\tilde{A}$, we are ready to define the dg category $\mathcal{C}(\tilde{A},\abimod)$ of dg $\tilde{A}$-bimodules in $\abimod$. An object $C=(C^q)_{q\in \zz}$ in $\mathcal{C}(\tilde{A},\abimod)$ is a graded $\tilde{A}$-bimodule in $\abimod$ endowed with a differential $d_C: C^q \rightarrow C^{q+1}$ of $A$-bimodules, satisfying the properties
\begin{equation} \label{axioms}
  d_C(f.x)=b_K(f).x + (-1)^p f.d_C(x), \ \ d_C(x.f)=d_C(x).f + (-1)^q x.b_K(f)
 \end{equation}
for any Koszul $p$-cochain $f$ with coefficients in $A$ and $x$ in $C^q$. The forgetful functor $\Und: \mathcal{C}(\tilde{A},\abimod) \rightarrow \mathcal{G}(\tilde{A},\abimod)$ is fully faithful.

By definition, for two objects $C$ and $C'$ in $\mathcal{C}(\tilde{A},\abimod)$, one has
$$\Hom_{\mathcal{C}(\tilde{A},\abimod)}(C,C')= \Hom_{\mathcal{G}(\tilde{A},\abimod)}(C,C'),$$
and this graded space has a differential $d$ defined by
\begin{equation} \label{diffofphi}
d(u)=d_{C'} \circ u - (-1)^r u \circ d_C 
\end{equation}
where $u : C \rightarrow C'$ is $r$-homogeneous. Then $(\Hom_{\mathcal{C}(\tilde{A},\abimod)}(C,C'),d)$ is a complex in $\fmod$.

The subcategory $\mathcal{C}_{str}(\tilde{A},\abimod)$ of $\mathcal{C}(\tilde{A},\abimod)$ has the same objects, with morphism spaces 
$$\Hom_{\mathcal{C}_{str}(\tilde{A},\abimod)}(C,C')= Z^0 (\Hom_{\mathcal{C}(\tilde{A},\abimod)}(C,C')),$$
that is, the morphisms from $C$ to $C'$ are of degree 0 and are morphisms of complexes. These morphisms are still called strict morphisms. The category $\mathcal{C}_{str}(\tilde{A},\abimod)$ is abelian.  The forgetful functor $\Und: \mathcal{C}_{str}(\tilde{A},\almod) \rightarrow \mathcal{G}_{str}(\tilde{A},\almod)$ is faithfully exact.

The homotopy category $\mathcal{K}(\tilde{A},\abimod)$ has the same objects as $\mathcal{C}(\tilde{A},\abimod)$, with morphism spaces
$$\Hom_{\mathcal{K}(\tilde{A},\abimod)}(C,C')= H^0 (\Hom_{\mathcal{C}(\tilde{A},\abimod)}(C,C')),$$
that is, the morphisms from $C$ to $C'$ are the homotopy classes of strict morphisms. Then $\mathcal{K}(\tilde{A},\abimod)$ is a triangulated category. 

From $\mathcal{K}(\tilde{A},\abimod)$, the derived category $\mathcal{D}(\tilde{A},\abimod)$ is constructed by inverting quasi-isomorphisms. One has the following canonical functors, the second one being triangulated.
$$P: \mathcal{C}_{str}(\tilde{A},\abimod) \longrightarrow \mathcal{K}(\tilde{A},\abimod),$$
$$Q: \mathcal{K}(\tilde{A},\abimod) \longrightarrow \mathcal{D}(\tilde{A},\abimod).$$

\Bdf \label{derkoscal2}
The data of the categories $\mathcal{C}(\tilde{A},\abimod)$ and $\mathcal{D}(\tilde{A},\abimod)$, and of the functor
$$\tilde{Q}=Q \circ P : \mathcal{C}_{str}(\tilde{A},\abimod) \longrightarrow \mathcal{D}(\tilde{A},\abimod)$$
is called the derived Koszul calculus of $A$. When the abelian category $\abimod$ is replaced by $\fmod$ in this definition, we get the scalar derived Koszul calculus of $A$.
\Edf

Let us remark that the scalar derived Koszul calculus of $A$ is the data of the categories $\mathcal{C}(\tilde{A})$ and $\mathcal{D}(\tilde{A})$, and of the corresponding functor
$$\tilde{Q}=Q \circ P : \mathcal{C}_{str}(\tilde{A}) \longrightarrow \mathcal{D}(\tilde{A}).$$
From the derived Koszul calculus to the scalar derived Koszul calculus, there are some obvious forgetful functors, such as $\mathcal{C}(\tilde{A},\abimod) \rightarrow \mathcal{C}(\tilde{A})$, obtained by forgetting the actions of $A$.

\subsection{The Koszul complex is an object in $\mathcal{C}(\tilde{A},\abimod)$} \label{K(A)dgbimod}

The most important example of dg $\tilde{A}$-bimodule in $\abimod$ is the Koszul complex $K(A)$~\cite[Proposition 3.6]{bt:kcpa}. The actions of a Koszul $p$-cochain $f: W_p \rightarrow A$ on an element $z=\alpha \otimes_k x_1 \ldots  x_q \otimes_k \beta$ in $K(A)_q$ are given by the following cap products
\begin{equation} \label{lcapkocomplex}
f \underset{K}{\frown} z = (-1)^{(q-p)p} \alpha \otimes_k x_1 \ldots  x_{q-p} \otimes_k f(x_{q-p+1} \ldots x_q)\beta,
\end{equation}
\begin{equation} \label{rcapkocomplex}
z \underset{K}{\frown} f = (-1)^{pq} \alpha f(x_1 \ldots x_p)\otimes_k x_{p+1} \ldots  x_q \otimes_k \beta
\end{equation}
when $p \leq q$. The right-hand sides are assumed to be zero when $p>q$.

The actions of $A$ on $z$ are given by $azb=a \alpha \otimes_k x_1 \ldots  x_q \otimes_k \beta b$ for any $a$ and $b$ in $A$. The actions of $\tilde{A}$ and $A$ on $K(A)$ commute since one has $f\underset{K}{\frown} (azb)=a(f \underset{K}{\frown} z)b$ and $(azb)\underset{K}{\frown} f=a(z \underset{K}{\frown} f)b$. Note that if $f$ is a $0$-cochain identified to an element $a \in \bigoplus_{i \in \mathcal{Q}_0} e_iAe_i$, then
$$f \underset{K}{\frown} z =  \alpha \otimes_k x_1 \ldots  x_q \otimes_k a \beta\ \ \mbox{and} \ \ az= a \alpha \otimes_k x_1 \ldots  x_q \otimes_k \beta$$
so that $f \underset{K}{\frown} z \neq az$ in many cases -- for example if $\mathcal{Q}$ has a single vertex and at least one arrow, choose $q=0$, $\alpha=\beta=1_A$ and $a \in \mathcal{Q}_1$.

Using the definition (\ref{definitiond}) of the differential $d$ of $K(A)$, it is direct to check that we have 
\begin{equation} \label{doncap1}
d(f \underset{K}{\frown} z) = b_K(f) \underset{K}{\frown} z + (-1)^p f \underset{K}{\frown} d(z),
\end{equation}
\begin{equation} \label{doncap2}
d(z \underset{K}{\frown} f) = d(z) \underset{K}{\frown} f + (-1)^q z \underset{K}{\frown} b_K(f).
\end{equation}

Another example of dg $\tilde{A}$-bimodule in $\abimod$ -- fundamental for defining an adapted Calabi-Yau property (Definition \ref{strongKcCY}) -- is given by $\Hom_{A^e}(K(A),A^e)$ defined in Proposition \ref{newobjects} just below.

\subsection{Some derived functors} \label{basicfunctors}

In Proposition \ref{fundamcomplexes}, we have recalled that $\Hom_{A^e}(K(A),M)$ and $M\otimes_{A^e} K(A)$ are dg $\tilde{A}$-bimodules. Let us construct more generally in Part (i) of the next proposition new objects in $\mathcal{C}(\tilde{A})$ from a given one in $\mathcal{C}(\tilde{A},\abimod)$. The construction in Part (ii) will be useful for the definition of $\Hom_{A^e}(K(A),A^e)$ as object in $\mathcal{C}(\tilde{A},\abimod)$. 

\Bpo \label{newobjects}
(i) For any chain (cochain) complex $C$ in $\mathcal{C}(\tilde{A},\abimod)$ and any $A$-bimodule $M$, $\Hom_{A^e}(C,M)$ (resp. $M\otimes_{A^e} C$) is a cochain complex in $\mathcal{C}(\tilde{A})$.

(ii) Choose in (i) the $A$-bimodule $M=A^e$ when $A^e$ is naturally considered as a left $A^e$-module. For this $A$-bimodule and for any chain complex $C$ in $\mathcal{C}(\tilde{A},\abimod)$, the cochain complex $\Hom_{A^e}(C,A^e)$ in $\mathcal{C}(\tilde{A})$ defined in (i) is underlying to a cochain complex in $\mathcal{C}(\tilde{A},\abimod)$. 
\Epo
\Bdm
Let us begin by the construction of $\Hom_{A^e}(C,M)$ in (i). Fix a chain complex $C=(C_q)_{q\in \zz}$ in $\mathcal{C}(\tilde{A},\abimod)$. Then $\Hom _{A^e} (C, M)$ is graded by the $\ff$-modules $\Hom _{A^e} (C, M)^q$ where $\Hom _{A^e} (C, M)^q:=\Hom _{A^e} (C_q, M)$. It is easy to check that $\Hom _{A^e} (C, M)$ is a graded $\tilde{A}$-bimodule for the following actions
$$(f.u)(x)=(-1)^p u(x.f), \ \ (u.f)(x)=u(f.x)$$
where $f:W_p \rightarrow A$, $u: C_q \rightarrow M$ and $x \in C_{q+p}$. Note that $x.f$ and $f.x$ are in $C_q$ by the graded actions of $\tilde{A}$ on $C$. In the first equality, $f$ of degree $p$ jumps over $u$ of degree $q$ and over $x$ of degree $q+p$, which gives the Koszul sign $(-1)^{p(q+q+p)}=(-1)^p$ -- while no Koszul sign occurs in the second equality.

Remark that if $C=K(A)$, we recover the cup actions, that is, $f.u=f\underset{K}{\smile}u$ and $u.f=u\underset{K}{\smile}f$.

Recall that if $d_C:C_{q+1}\rightarrow C_q$ is the differential of $C$, the differential
$$\Hom_{A^e}(d_C,M): \Hom _{A^e} (C_q,M) \rightarrow \Hom _{A^e} (C_{q+1},M)$$
of $\Hom _{A^e} (C,M)$ is defined by $\Hom_{A^e}(d_C,M)(u)=-(-1)^q u\circ d_C$ for $u \in \Hom _{A^e} (C_q,M)$. The reader can carefully verify the following
\begin{equation} \label{verifaxiom1}
\Hom_{A^e}(d_C,M)(f.u)=b_K(f).u + (-1)^p f.\Hom_{A^e}(d_C,M)(u), 
 \end{equation}
\begin{equation} \label{verifaxiom2}
\Hom_{A^e}(d_C,M)(u.f)=\Hom_{A^e}(d_C,M)(u).f + (-1)^q u.b_K(f),
 \end{equation}
so that $\Hom_{A^e}(C,M)$ is a dg $\tilde{A}$-bimodule.

Similarly, for any cochain dg $\tilde{A}$-bimodule $C$ in $\abimod$, $M\otimes_{A^e} C$ is a cochain dg $\tilde{A}$-bimodule for the following actions
$$f.(m \otimes_{A^e} x)=m \otimes_{A^e} (f.x), \ \ (m \otimes_{A^e} x).f=m \otimes_{A^e} (x.f)$$
where $f \in \tilde{A}$, $m\in M$ and $x \in C$. The definition of the differential of $M\otimes_{A^e} C$ from the differential of $C$ is clear and the verification of the properties is direct.  

(ii) Now $A^e$ is the left $A^e$-module defined by the product of $A^e$, obtaining an $A$-bimodule $M=A^e$. Let us begin to define on the dg $\tilde{A}$-bimodule $\Hom_{A^e}(C,A^e)$ constructed above a structure of $A$-bimodule. Fix $u \in \Hom _{A^e} (C_q,A^e)$. Note that $u$ is $A^e$-linear by assumption, meaning that for any $a$, $b$ in $A$
and $x\in C_q$, one has
\begin{equation} \label{leftA^elinear}
u(a x b)= (a \otimes b) . u(x),
\end{equation}
where the dot denotes the product in the algebra $A^e$.

For $\alpha$ and $\beta$ in $A$, define $\alpha u \beta :C_q \rightarrow A^e$ by
\begin{equation} \label{rightA^elinear}
(\alpha u \beta )(x)= u(x) . (\beta \otimes \alpha)
\end{equation} 
for any $x\in C_q$, with the same meaning for the dot. Then
$$(\alpha u \beta) (a x b) = ((a\otimes b).u(x)). (\beta \otimes \alpha) = (a\otimes b).(u(x). (\beta \otimes \alpha)) = (a\otimes b).((\alpha u \beta )(x))$$
since the product in $A^e$ is associative. So $\alpha u \beta$ is in $\Hom _{A^e} (C_q,A^e)$.

Next it is easy to check that
$$\alpha' (\alpha u \beta ) \beta'= (\alpha' \alpha) u (\beta \beta')$$
where $\alpha'$ and $\beta'$ are in $A$, as well that $1_A u 1_A=u$. So $\Hom_{A^e}(C,A^e)$ is an $A$-bimodule.

Recall that for $f:W_p \rightarrow A$, $u: C_q \rightarrow A^e$ and $x \in C_{q+p}$
, one has 
$$(f.u)(x)=(-1)^p u(x.f), \ \ (u.f)(x)=u(f.x).$$
It is the immediate to check that
$$f.(\alpha u \beta )=\alpha (f.u) \beta \ \ \mathrm{and} \ \ (\alpha u \beta ).f=\alpha (u.f) \beta,$$
so that $\Hom_{A^e}(C,A^e)$ is a graded $\tilde{A}$-bimodule in $\abimod$.

In order to show that $\Hom_{A^e}(C,A^e)$ is an object in $\mathcal{C}(\tilde{A},\abimod)$, it remains to prove that the differential $\Hom_{A^e}(d_C,A^e)$ is a morphism of $A$-bimodule. From the definition
$\Hom_{A^e}(d_C,A^e)(\alpha u \beta)=-(-1)^q (\alpha u \beta) \circ d_C$
applied to any $x \in C_{q+1}$ and using the definition (\ref{rightA^elinear}) of $\alpha u \beta$, we arrive easily to
$$\Hom_{A^e}(d_C,A^e)(\alpha u \beta)= \alpha \Hom_{A^e}(d_C,A^e)(u) \beta. \ \Edm$$

As a consequence of Proposition \ref{newobjects}, for any $A$-bimodule $M$, we define the following covariant functors
\begin{align} \label{firstfunct}
\Hom_{A^e}(-,M) :\ \  & \mathcal{C}(\tilde{A},\abimod)^{op} \longrightarrow \mathcal{C}(\tilde{A}) \\ \label{secondfunct}
M\otimes_{A^e} - :\  \ & \mathcal{C}(\tilde{A},\abimod) \longrightarrow \mathcal{C}(\tilde{A})\\ \label{thirdfunct}
\Hom_{A^e}(-,A^e) :\ \  & \mathcal{C}(\tilde{A},\abimod)^{op} \longrightarrow \mathcal{C}(\tilde{A},\abimod)
\end{align}
where $\mathcal{C}(\tilde{A},\abimod)^{op}$ denotes the opposite dg category~\cite[Definition 3.9.2]{yeku:dercat}. The definition of these functors on morphisms and the verification of the functoriality propertiess are direct. Let us make precise that for the first functor. It is similar for the two other functors.

Let $C$ and $C'$ be chain dg $\tilde{A}$-bimodules in $\abimod$. Let $\phi : C \rightarrow C'$ be of degree $r$ in $\mathcal{C}(\tilde{A},\abimod)$, that is, $\phi$ is an $r$-homogeneous graded morphism of $\tilde{A}$-bimodules and each $\phi: C_q \rightarrow C'_{q+r}$ is a morphism of $A$-bimodules. We define the graded linear map
$$\Hom_{A^e}(\phi,M): \Hom_{A^e}(C',M) \longrightarrow \Hom_{A^e}(C,M)$$
of degree $-r$ by
\begin{equation} \label{functoronmorphism}
\Hom_{A^e}(\phi,M)(u)=(-1)^{r(q+r)} u\circ \phi \in \Hom_{A^e}(C_q,M)
\end{equation}
for $u:C'_{q+r} \rightarrow M$, since $\phi: C_q \rightarrow C'_{q+r}$. The sign is a Koszul sign.

\Brm \label{importantrem}
Equality (\ref{functoronmorphism}) extended by the same formula to $d_C:C \rightarrow C$ of degree $-1$ (only morphism of $A$-bimodules) coincides with the opposite of the differential $\Hom_{A^e}(d_C,M)$ defined and used in the proof of Proposition \ref{newobjects}.
\Erm

Then it is direct to check that, if $f:W_p \rightarrow A$ is a Koszul $p$-cochain, one has
\begin{align}
\Hom_{A^e}(\phi,M)(f.u) & = (-1)^{rp} f.\Hom_{A^e}(\phi,M)(u)\\
\Hom_{A^e}(\phi,M)(u.f) & = \Hom_{A^e}(\phi,M)(u).f
\end{align}
where we get a Koszul sign in the first equality, and none in the second equality. So $\Hom_{A^e}(\phi,M)$ is a morphism in $\mathcal{C}(\tilde{A})$ and ($\ref{firstfunct}$) is defined on morphisms.

The first functoriality property $\Hom_{A^e}(\id_C,M)=\id_{\Hom_{A^e}(C,M)}$ is immediate. It is direct to verify the second one
\begin{equation} \label{functoriality2}
  \Hom_{A^e}(\psi \circ \phi,M)=(-1)^{rr'}\Hom_{A^e}(\phi,M)\circ \Hom_{A^e}(\psi,M)
\end{equation}
where $\phi :C\rightarrow C'$ has degree $r$ and $\psi :C' \rightarrow C''$ has degree $r'$. In conclusion, we have defined the functor ($\ref{firstfunct}$).

\Bpo \label{functareDG}
The three functors $(\ref{firstfunct})$-$(\ref{thirdfunct})$ are dg functors.
\Epo
\Bdm
Let us limit ourselves to the contravariant functor $F=\Hom_{A^e}(-,M)$. The proof for the functors $(\ref{secondfunct})$ and $(\ref{thirdfunct})$ is left to the reader. Let $\phi:C \rightarrow C'$ be a morphism of degree $r$ in $\mathcal{C}(\tilde{A},\abimod)$, where the differentials of the chain complexes $C=(C_q)_{q\in \zz}$ and $C'=(C'_q)_{q\in \zz}$ are denoted by $d_C$ and $d_{C'}$ respectively. Recall that the graded space $\Hom_{\mathcal{C}(\tilde{A},\abimod)}(C,C')$ has a differential $d$ defined by
\begin{equation} \label{diffofphi2}
d(\phi)=d_{C'} \circ \phi - (-1)^r \phi \circ d_C.
\end{equation}
Similarly, $\Hom_{\mathcal{C}(\tilde{A})}(F(C'),F(C))$ has a differential $d'$ defined by
\begin{equation} \label{diffofpsi}
d'(\psi)= \Hom_{A^e}(d_C,M) \circ \psi - (-1)^r \psi \circ \Hom_{A^e}(d_{C'},M),
\end{equation}
where $\psi: F(C') \rightarrow F(C)$ has degree $r$ in $\mathcal{C}(\tilde{A})$. Choose $\psi=F(\phi)$ which has degree $-r$. According to~\cite[Definition 3.9.1]{yeku:dercat}, $F$ is a dg functor if one has
\begin{equation} \label{verifDGfunctor}
F(d(\phi))=d'(F(\phi)).
\end{equation}
To prove that, apply the functoriality property (\ref{functoriality2}) to (\ref{diffofphi2}) and get
$$F(d(\phi))=(-1)^r F(\phi)\circ F(d_{C'}) - (-1)^r (-1)^r F(d_C) \circ F(\phi).$$
Then using Remark \ref{importantrem}, one has $F(d_C)=-\Hom_{A^e}(d_C,M)$ and $F(d_{C'})=-\Hom_{A^e}(d_{C'},M)$. Therefore $F(d(\phi))= d'(F(\phi))$ by (\ref{diffofpsi}).
\Edm

\Bcr \label{passtotriang}
The three dg functors $(\ref{firstfunct})$-$(\ref{thirdfunct})$ pass to homotopy categories, providing triangulated functors.
\Ecr
  
It is an immediate consequence of~\cite[Theorem 5.6.10, Theorem 5.5.1]{yeku:dercat}. So we obtain the triangulated functors
\begin{align} \label{threetriangfunc1}
  \Hom_{A^e}(-,M) :\ \  & \mathcal{K}(\tilde{A},\abimod)^{op} \longrightarrow \mathcal{K}(\tilde{A}) \\
\label{threetriangfunc2}
M\otimes_{A^e} - :\  \ & \mathcal{K}(\tilde{A},\abimod) \longrightarrow \mathcal{K}(\tilde{A}) \\
\label{threetriangfunc3}
\Hom_{A^e}(-,A^e) :\ \  & \mathcal{K}(\tilde{A},\abimod)^{op} \longrightarrow \mathcal{K}(\tilde{A},\abimod).
\end{align}

Postcomposing by a localization functor $Q$, we get the triangulated functors
\begin{align} \label{threeQtriangfunc1}
  \Hom_{A^e}(-,M) :\ \  & \mathcal{K}(\tilde{A},\abimod)^{op} \longrightarrow \mathcal{D}(\tilde{A}) \\
\label{threeQtriangfunc2}
M\otimes_{A^e} - :\  \ & \mathcal{K}(\tilde{A},\abimod) \longrightarrow \mathcal{D}(\tilde{A}) \\
\label{threeQtriangfunc3}
\Hom_{A^e}(-,A^e) :\ \  & \mathcal{K}(\tilde{A},\abimod)^{op} \longrightarrow \mathcal{D}(\tilde{A},\abimod).
\end{align}

By Proposition $\ref{fundaresult}$, the category $\mathcal{C}(\tilde{A}, \abimod)$ has enough $K$-projectives. Therefore, Theorem $\ref{contrarderfunctexist}$ and Theorem $\ref{lderfunctexist}$ show that (\ref{threeQtriangfunc1})-(\ref{threeQtriangfunc3}) have triangulated derived functors 
\begin{align} \label{threederfunc1}
  \RHom_{A^e}(-,M) :\ \  & \mathcal{D}(\tilde{A},\abimod)^{op} \longrightarrow \mathcal{D}(\tilde{A}) \\
\label{threederfunc2}
M\stackrel{L}{\otimes}_{A^e} - :\  \ & \mathcal{D}(\tilde{A},\abimod) \longrightarrow \mathcal{D}(\tilde{A}) \\
\label{threederfunc3}
\RHom_{A^e}(-,A^e) :\ \  & \mathcal{D}(\tilde{A},\abimod)^{op} \longrightarrow \mathcal{D}(\tilde{A},\abimod).
\end{align}
Note that for the existence of (\ref{threederfunc2}), it suffices that the category $\mathcal{C}(\tilde{A}, \abimod)$ has enough K-flat objects. But we know that any K-projective object is K-flat~\cite[Proposition 10.3.4]{yeku:dercat}.

\subsection{Defining strong Kc-Calabi-Yau algebras} \label{definingstrongKcCY}

Recall that a quadratic quiver algebra $A$ over a finite quiver $\mathcal{Q}$ is said to be \emph{$n$-Kc-Calabi-Yau} if the $A$-bimodule Koszul complex $K(A)$ of $A$ has finite length $n$ and if 
$$\RHom_{A^e}(K(A), A^e) \cong K(A)[-n]$$
in the bounded derived category $\mathcal{D}^b(\abimod)$ of $A$-bimodules~\cite[Definition 5.1]{bt:kcpa}. The abbreviation $Kc$ stands for \emph{Koszul complex}. 

If $A$ is $n$-Kc-Calabi-Yau and Koszul, then $A$ is homologically smooth since $K(A)$ is a finitely generated projective bimodule resolution of $A$ having a finite length. Moreover $K(A) \cong A$ in $\mathcal{D}^b(\abimod)$, so that one has
$$\RHom_{A^e}(A, A^e) \cong A[-n] \ \mathrm{in} \ \mathcal{D}^b(\abimod),$$
meaning that the associative algebra $A$ is $n$-Calabi-Yau in Ginzburg's sense~\cite{vg:cy}.

For an $n$-Kc-Calabi-Yau algebra $A$, the spaces $\HK^p(A,M)$ and $\HK_{n-p}(A,M)$ are isomorphic for any $A$-bimodule $M$ and $0\leq p \leq n$~\cite[Theorem 5.5]{bt:kcpa}. When $A$ is moreover Koszul, we recover the Van den Bergh duality theorem expressed in terms of Hochschild (co)homology and particularized to $n$-Calabi-Yau quadratic algebras~\cite{vdb:dual}.

In~\cite{bt:kcpa}, it is proved that the preprojective algebra of a connected graph distinct from the types $\mathrm{A}_1$ and $\mathrm{A}_2$ is 2-Kc-Calabi-Yau. Such a preprojective algebra is known to be \emph{not 2-Calabi-Yau} in Ginzburg's sense when the graph is Dynkin ADE, because the minimal projective resolution of $A$ has then an infinite length so that $A$ is not homologically smooth. If the graph is not Dynkin ADE, the preprojective algebra is Koszul.

Since the derived functor ($\ref{threederfunc3}$) exists, the derived functor
$$\RHom_{A^e}(-,A^e) : \mathcal{D}^b(\tilde{A},\abimod)^{op} \rightarrow \mathcal{D}^b(\tilde{A},\abimod)$$
exists. Here the bounded derived category $\mathcal{D}^b(\tilde{A},\abimod)$ is defined by inverting quasi-isomorphisms from the full triangulated subcategory $\mathcal{K}^b(\tilde{A},\abimod)$ of $\mathcal{K}(\tilde{A},\abimod)$ consisting of \emph{bounded} graded objects~\cite[Definition 7.3.3]{yeku:dercat}. Note that the functor $\Hom_{A^e}(-,A^e)$ sends bounded graded objects to bounded graded objects.

\Brm \label{twoboundedness}
If $K(A)$ has a finite length, i.e., if the chain complex $K(A)$ is bounded below, the grading of the Koszul cochain dg algebra $\tilde{A}$ is bounded above. Then $\mathcal{D}^b(\tilde{A},\abimod)$ is canonically equivalent to the full subcategory of $\mathcal{D}(\tilde{A},\abimod)$ on the objects $C$ whose cohomology $H(C)$ is bounded~\cite[Definition 7.3.4]{yeku:dercat}. In fact, it is direct to extend~\cite[Proposition 7.3.12]{yeku:dercat} from nonpositive dg algebras to bounded above dg algebras. Same remark for $\fmod$ instead of $\abimod$.
\Erm

We are now ready to introduce the following definition. Actually, this definition was given in~\cite[Definition 5.10]{bt:kcpa} with the further assumption that the contravariant endofunctor $\RHom_{A^e}(-,A^e)$ of $\mathcal{D}^b(\tilde{A},\abimod)$ exists, but this existence assumption is unnecessary.  

\Bdf \label{strongKcCY}
Let $A$ be a quadratic quiver algebra over a finite quiver $\mathcal{Q}$. We say that $A$ is strong $n$-Kc-Calabi-Yau if the complex $K(A)$ has finite length $n$ and if
$$\RHom_{A^e}(K(A), A^e) \cong K(A)[-n]$$
in the bounded derived category $\mathcal{D}^b(\tilde{A},\abimod)$.
\Edf

If $A$ is strong $n$-Kc-Calabi-Yau, then $A$ is $n$-Kc-Calabi-Yau by forgetting the actions of $\tilde{A}$. The above definition was motivated by the following result.

\Bpo \label{preproarestrong}
Let $A$ be the preprojective algebra of a connected graph distinct from the types $\mathrm{A}_1$ and $\mathrm{A}_2$. Then $A$ is strong $2$-Kc-Calabi-Yau.
\Epo

The proof is the same as in~\cite[Subsection 5.3]{bt:kcpa} since the required existence assumption is now unnecessary. We prove in the next Section that a polynomial algebra in $n$ variables is strong $n$-Kc-Calabi-Yau.

\subsection{A strong Poincar\'e Van den Bergh duality theorem} \label{strongduality}

If $A$ is $n$-Kc-Calabi-Yau, the graded spaces $\HK^{\bullet}(A)$ and $\HK_{n-\bullet}(A)$ are isomorphic. In the strong case, such an isomorphism is expressed as a cap action by a fundamental class like the genuine Poincar\'e duality~\cite[Theorem 3.30]{hat:algtop}.

\Bte \label{strongvdbduality2}
Let $A$ be a quadratic quiver algebra over a finite quiver $\mathcal{Q}$. Assume that $A$ is strong $n$-Kc-Calabi-Yau. Then there is a class $c$ in $\HK_n(A)$, called fundamental class, such that
$$c \underset{K}{\frown} -  : \HK^{\bullet}(A) \rightarrow \HK_{n-\bullet}(A)$$
is an isomorphism of $\HK^{\bullet}(A)$-bimodules, inducing an isomorphism of $\HK^{\bullet}_{hi}(A)$-bimodules from $\HK^{\bullet}_{hi}(A)$ to $\HK^{hi}_{n-\bullet}(A)$. For all $\alpha \in \HK^p(A)$, one has $c \underset{K}{\frown} \alpha = (-1)^{np} \alpha \underset{K}{\frown}c$.
\Ete

The proof is the same as the proof given for~\cite[Theorem 5.11]{bt:kcpa}, with the same definition of the fundamental class $c$ and of the higher Koszul calculus $\HK^{\bullet}_{hi}(A)$ and $\HK^{hi}_{\bullet}(A)$. However, in the statement of~\cite[Theorem 5.11]{bt:kcpa}, it is assumed that the derived functors 
\begin{align} \label{neededfunctors}
\RHom_{A^e}(-,A) :\ \  & \mathcal{D}^b(\tilde{A},\abimod)^{op} \longrightarrow \mathcal{D}^b(\tilde{A}) \\
A\stackrel{L}{\otimes}_{A^e} - :\  \ & \mathcal{D}^b(\tilde{A},\abimod) \longrightarrow \mathcal{D}^b(\tilde{A})
\end{align}
exist. Now these derived functors exist because functors ($\ref{threederfunc1}$) and ($\ref{threederfunc2}$) exist for $M=A$. So the proof is complete.

If $A$ is the preprojective algebra of a connected graph distinct from the types $\mathrm{A}_1$ and $\mathrm{A}_2$, the fundamental class $c$ is explicitly known in terms of the graph~\cite[Subsection 4.3]{bt:kcpa}. The conclusion of Theorem \ref{strongvdbduality2} was proved directly in this case~\cite[Theorem 4.4]{bt:kcpa}, getting a duality theorem before stating a general theorem. If the graph is moreover Dynkin ADE, the algebras $\HK^{\bullet}(A)$ and $\HK^{\bullet}_{hi}(A)$ were explicitly computed by generators and relations, so that the duality theorem allowed to explicitly compute the $\HK^{\bullet}(A)$-bimodule $\HK_{\bullet}(A)$. See~\cite[Section 6]{bt:kcpa} for the result of these computations. These computations show that the Koszul calculus provides more cohomological invariants than the Hochschild calculus, except in type $\mathrm{E}_8$ with $\car(\ff)=2$ where the invariants are the same.

\setcounter{equation}{0}

\section{A polynomial algebra in $n$ variables is strong $n$-Kc-Calabi-Yau} \label{polynomials}

The goal of this section is to obtain a duality theorem for polynomial algebras implying that these algebras are strong Kc-Calabi-Yau. This duality theorem is similar to the duality theorem obtained for preprojective algebras~\cite[Theorem 4.4]{bt:kcpa}.

\subsection{Recalling some properties of a symmetric algebra} \label{symmalg}

In this Section, we apply the setup of Subsection \ref{koscalA} to a quiver $\mathcal{Q}$ having a single vertex and exactly $n$ arrows with $n\geq 1$. So $k=\ff$ and $\dim V =n$. In the sequel, the base field will be denoted by $k$ instead of $\ff$. The characteristic of $k$ is arbitrary.

We are interested in the \emph{symmetric algebra} $A=S(V)$ of the vector space $V$. By definition, $A$ is the quadratic algebra $A=T(V)/(R)$, where $R$ is the subspace of $V \otimes V$ generated by the elements $x\otimes y -y \otimes x$ for $x$ and $y$ in $V$. We can see $A$ as a polynomial algebra in $n$ variables (the arrows of the quiver) over $k$, but it is better to express the duality theorem intrinsically, that is, in terms of the symmetric algebra $A=S(V)$. It is classical that $A$ is Koszul~\cite{weib:homo} and $n$-Calabi-Yau~\cite{vg:cy}, therefore $A$ is $n$-Kc-Calabi-Yau.

In our situation, let us give an explicit description of the subspace $W_p$ of $V^{\otimes p}$ defined as in Equality (\ref{definitionw}) by 
$$W_{p}=\bigcap_{i+2+j=p}V^{\otimes i}\otimes R\otimes V^{\otimes j}.$$

We need the exterior algebra $\Lambda(V)$ of $V$. It is the quadratic algebra defined by $\Lambda(V)=T(V)/(R')$ where $R'$ is the subspace of $V \otimes V$ generated by the elements $x\otimes x$ for $x$ in $V$. The algebra $\Lambda(V)$ is graded by the subspaces $\Lambda^p(V)$, $p\geq 0$. The following description is standard~\cite[A III.190]{bou:alg} and is valid whether $V$ is finite-dimensional or not.

\Bpo \label{standardbou}
Let $V$ be a $k$-vector space and $A=S(V)$ be the symmetric algebra of $V$. For any $p\geq 0$, the space $W_p$ is equal to the image of the $k$-linear map
$Ant_p:V^{\otimes p}\rightarrow V^{\otimes p}$ defined by
$$Ant_p(v_1 \otimes \ldots \otimes v_p)= \sum _{\sigma \in \Sigma_p} \mathrm{sgn}(\sigma)\, v_{\sigma^{-1} (1)} \otimes \ldots \otimes v_{\sigma^{-1} (p)}$$
for any $v_1, \ldots, v_p$ in $V$, where $\Sigma_p$ is the symmetric group of $\{1, ... , p\}$ and $\mathrm{sgn}(\sigma)$ is the sign of $\sigma$.
Moreover, $Ant_p$ induces a linear isomorphism $\Lambda^p (V) \rightarrow W_p$. When $V$ is $n$-dimensional, $W_p=0$ for $p>n$.
\Epo

An $A$-bimodule $M$ is said to be $A$-central if left action and right action of $A$ over $M$ are equal -- that makes sense since the algebra $A$ is commutative. For example the $A$-bimodule $A$ is $A$-central, but $A\otimes A$ is not $A$-central. The following is the same as Proposition 3.16 in~\cite{bls:kocal} and is again valid whether $V$ is finite-dimensional or not.

\Bpo \label{symalg}
Let $V$ be a $k$-vector space and $A=S(V)$ be the symmetric algebra of $V$. Let $M$ be an $A$-central $A$-bimodule. 
The differentials $b^K$ and $b_K$ of the complexes $M\otimes W_{\bullet}$ and $\Hom(W_{\bullet},M)$ are zero. Therefore, $\HK_{\bullet}(A,M)=M\otimes W_{\bullet}$ and $\HK^{\bullet}(A,M)=\Hom(W_{\bullet},M)$. In particular, for $M=A$, the Koszul cochain dg algebra $\tilde{A}$ (Definition \ref{defAtilde}) has a zero differential.
\Epo
\Bdm
We only prove that $b_K=0$. For a proof of $b^K =0$, see~\cite[Proposition 3.16]{bls:kocal}. For $f\in \Hom(W_p,M)$, $b_K(f)\in \Hom(W_{p+1},M)$ is defined by the Equality (\ref{defb_K}), that is, by
\begin{equation} \label{defb_K2}
  b_K(f)( x_1 \ldots x_{p+1}) =f(x_1\ldots x_{p})\, x_{p+1} -(-1)^p x_1\, f(x_2 \ldots x_{p+1})
\end{equation}
where $x_1 \ldots x_{p+1} \in W_{p+1}$. Recall that the notation $x_1 \ldots x_{p+1}$ in the left-hand side is conventional for abbreviating a sum of such elements. Moreover in the right-hand side, we keep the same convention, so that $x_1\ldots x_{p}$ in the first term and $x_2 \ldots x_{p+1}$ in the second term are viewed in $W_p$.

Using Proposition \ref{standardbou}, we have to compute $b_K(f)(Ant_{p+1}(v_1 \ldots  v_{p+1}))$ for any $v_1, \ldots, v_{p+1}$ in $V$. In $Ant_{p+1}(v_1 \ldots  v_{p+1})$, the symbols $\otimes$ are omitted for brevity. The standard formula
\begin{equation}  \label{standardformula1}
Ant_{p+1}(v_1  \ldots  v_{p+1}) =\sum_{1 \leq i \leq p+1} (-1)^{p+1-i} Ant_p(v_1  \ldots  \widehat{v_i} \dots  v_{p+1}) \otimes v_i
\end{equation}
encodes the inclusion $W_{p+1} \subseteq W_p \otimes V$, while the standard formula 
\begin{equation}  \label{standardformula2}
Ant_{p+1}(v_1  \ldots  v_{p+1}) =\sum_{1\leq i \leq p+1} (-1)^{i-1} v_i \otimes Ant_p(v_1  \ldots  \widehat{v_i} \dots  v_{p+1})
\end{equation}
encodes the inclusion $W_{p+1} \subseteq V \otimes W_p$. Here $\widehat{v_i}$ means that $v_i$ is omitted. Therefore, applying (\ref{defb_K2}), we obtain
\begin{align} 
  b_K(f)(Ant_{p+1}(v_1 \ldots  v_{p+1})) = & \sum_{1 \leq i \leq p+1} (-1)^{p+1-i} \ [ \, f(Ant_p(v_1  \ldots  \widehat{v_i} \dots  v_{p+1})) v_i  \nonumber \\
  & - v_i f(Ant_p(v_1  \ldots  \widehat{v_i} \dots  v_{p+1}))\, ]. \nonumber
\end{align}
Thus $b_K(f)(Ant_{p+1}(v_1 \ldots  v_{p+1}))=0$ since $M$ is $A$-central.
\Edm
\\

Since $A$ is Koszul, Proposition \ref{symalg} provides an explicit description of $\HH_{\bullet}(A,M)$ and $\HH^{\bullet}(A,M)$. When $\dim V=n$ and $M=A$, we obtain as in~\cite{weib:homo} well-known graded A-module isomorphims  
$$\HH_{\bullet}(A) \cong \Lambda_A^{\bullet}(A^n) \cong \HH^{\bullet}(A).$$

\subsection{A duality isomorphism for $A=S(V)$} \label{dualitysymmalg}

Let us keep the assumptions of the previous subsection, that is, $A=S(V)$ with $\dim V=n$. Let us fix a basis $\{x_1, \ldots, x_n\}$ of the space $V$. According to Proposition \ref{standardbou}, for $0\leq p \leq n$, the space $W_p$ is identified to $\Lambda^p (V)$ and the following elements
\begin{equation}  \label{W_pbasis}
x_{i_1} \wedge \ldots  \wedge x_{i_p}, \ 1\leq i_1 < \ldots i_p \leq n,
\end{equation}
form a basis of $W_p$. In particular, $W_n$ is one-dimensional generated by $x_1 \wedge \ldots  \wedge x_n$. Proposition \ref{symalg} shows that $\HK_n(A)=A\otimes W_n$ so that $\HK_n(A)$ is a free $A$-module generated by $1_A \otimes (x_1 \wedge \ldots \wedge x_n)$.

\Bdf \label{fundclass}
The element $c=1_A \otimes (x_1 \wedge \ldots \wedge x_n)$ in $\HK_n(A)$ is called a \emph{fundamental class} of $A$.
\Edf

In Theorem \ref{dualityiso} below, the duality isomorphism is valid for any $A$-bimodule $M$, $A$-central or not, so that the differentials of the complexes $M\otimes W_{\bullet}$ and $\Hom(W_{\bullet},M)$ are not necessarily zero. Let us begin by a lemma.

\Blm 
For all Koszul $p$-cochain $f:W_p \rightarrow M$, one has
\begin{equation} \label{comformula}
  c \underset{K}{\frown} f= (-1)^{np} f \underset{K}{\frown} c.
\end{equation}
\Elm
\Bdm
Apply the general definition of the Koszul cap product~\cite{bt:kcpa} to $c=1_A \otimes (y_1\ldots y_n)$ in $A \otimes W_n$ and to the Koszul $p$-cochain $f$. We obtain a Koszul $(n-p)$-chain $c \underset{K}{\frown} f$ in $M \otimes W_{n-p}$ given by
$$c \underset{K}{\frown} f =(-1)^{np} f(y_1 \ldots y_p)\otimes (y_{p+1} \ldots y_n)$$
following the conventional notation for $y_1\ldots y_n$ as in Formula (\ref{defb_K2}). Via the identifications $W_n=\Lambda^n(V)$ and $W_{n-p}=\Lambda^{n-p}(V)$ and using the basis $\{x_1, \ldots, x_n\}$ of $V$, we get
\begin{eqnarray} \label{rightcap1}
c\underset{K}{\frown} f = (-1)^{np} \sum _{\sigma \in Sh(p,n-p)} \mathrm{sgn}(\sigma) \,f(x_{\sigma(1)} \wedge \ldots \wedge x_{\sigma(p)})\nonumber \\
\otimes \, (x_{\sigma(p+1)}\wedge \ldots \wedge x_{\sigma(n)}),   
\end{eqnarray}
where $Sh(p,n-p)$ is the set of the $(p,n-p)$-shuffles, i.e., permutations $\sigma$ of $\{1,\ldots ,n\}$ such that $\sigma(1)< \cdots < \sigma(p)$ and $\sigma(p+1)< \cdots < \sigma(n)$.

Similarly, the general definition
$$f \underset{K}{\frown} c =(-1)^{(n-p)p} f(y_{n-p+1} \ldots y_n)\otimes (y_1 \ldots y_{n-p})$$
gives via similar identifications
\begin{eqnarray} \label{leftcap1}
f\underset{K}{\frown} c = (-1)^{(n-p)p} \sum _{\sigma \in Sh(n-p,p)} \mathrm{sgn}(\sigma) \,f(x_{\sigma(n-p+1)} \wedge \ldots \wedge x_{\sigma(n)})\nonumber \\
\otimes \, (x_{\sigma(1)}\wedge \ldots \wedge x_{\sigma(n-p)}),  
\end{eqnarray}
where $Sh(n-p,p)$ is the set of the $(n-p,p)$-shuffles. Then it is direct to deduce from (\ref{leftcap1}) the following
\begin{eqnarray} \label{leftcap2}
f\underset{K}{\frown} c = \sum _{\sigma \in Sh(p,n-p)} \mathrm{sgn}(\sigma) \,f(x_{\sigma(1)} \wedge \ldots \wedge x_{\sigma(p)})\nonumber \\
\otimes \, (x_{\sigma(p+1)}\wedge \ldots \wedge x_{\sigma(n)}).  
\end{eqnarray}
Comparing with (\ref{rightcap1}), we obtain (\ref{comformula}).
\Edm
\\

Note that the expression of $f\underset{K}{\frown} c$ in (\ref{leftcap2}) has a classical meaning. In fact, when $M=A$, $f \underset{K}{\frown} c=i_f(x_1 \wedge \ldots \wedge x_n)$ where $i_f$ is the contraction map (inner or internal product) used in differential geometry~\cite{lpv:poisson}.

\Bte \label{dualityiso}
Let $V$ be a $k$-vector space of finite dimension $n$ and $A=S(V)$ be the symmetric algebra of $V$. Let $M$ be an $A$-bimodule. For each Koszul $p$-cochain $f$ with coefficients in $M$, we define the Koszul $(n-p)$-chain $\theta_M (f)$ with coefficients in $M$ by
\begin{equation} \label{deftheta}
\theta_M (f) = c \underset{K}{\frown} f.
\end{equation}
Then the linear map $\theta_M: \Hom(W_{\bullet},M) \rightarrow M\otimes W_{n-\bullet}$ is a $(-n)$-degree complex isomorphism of dg $\tilde{A}$-bimodules.

In particular, $\theta_A$ realizes a graded algebra isomorphism from $\tilde{A}$ to the tensor product $A \otimes \Lambda(W_{n-1})$ of the algebra $A$ and the exterior algebra $\Lambda(W_{n-1})$ of the $n$-dimensional space $W_{n-1}$ considered in degree $1$.
\Ete
\Bdm
The linear map $\theta_M$ is graded homogeneous of degree $-n$ since it sends a $p$-cochain $f$ to a $(p-n)$-cochain.

Let $N$ be an $A$-bimodule and $g$ a Koszul cochain with coefficients in $N$. Using general \emph{associativity relations} between Koszul cup and cap products~\cite{bt:kcpa} and Equality (\ref{comformula}), one has
$$\theta_{M\otimes_A N} (f\underset{K}{\smile} g)= c \underset{K}{\frown} (f\underset{K}{\smile} g)= (c \underset{K}{\frown} f) \underset{K}{\frown} g=(-1)^{np} (f \underset{K}{\frown} c ) \underset{K}{\frown} g=(-1)^{np}f \underset{K}{\frown} (c  \underset{K}{\frown} g).$$

In particular, for $N=A$ and for $M=A$ successively, one has
\begin{equation} \label{theta_M}
\theta_{M} (f\underset{K}{\smile} g) = \theta_M (f) \underset{K}{\frown} g,
\end{equation}
\begin{equation}  \label{theta_N}
\theta_{N} (f\underset{K}{\smile} g) =  (-1)^{np} f \underset{K}{\frown} \theta_N (g).
\end{equation}
Exchanging $N$ and $M$, $f$ and $g$ in the latter equality, one has
\begin{equation}  \label{exchangetheta_N}
\theta_{M} (g\underset{K}{\smile} f) =  (-1)^{nq} g \underset{K}{\frown} \theta_M (f).
\end{equation}
where $f: W_{\bullet} \rightarrow M$ and $g: W_q \rightarrow A$. Then (\ref{theta_M}) and (\ref{exchangetheta_N}) show that $\theta_M: \Hom(W_{\bullet},M) \rightarrow M\otimes W_{n-\bullet}$ is a $(-n)$-degree morphism of graded $\tilde{A}$-bimodules, where $\tilde{A}=\Hom(W_{\bullet},A)$ is just considered as a graded algebra for the Koszul cup product $\underset{K}{\smile}$.

In order to get a complex morphism of dg $\tilde{A}$-bimodules, it remains to examine what happens for the Koszul differentials. Equation (\ref{theta_N}), together with  (\ref{theta_M}) after exchanges, implies that
\begin{equation} \label{thetaoncommutator}
\theta_N ([f,g]_{\underset{K}{\smile}}) =  (-1)^{np} [f,\theta_N (g)]_{\underset{K}{\frown}}
\end{equation} 
for any $f:W_p \rightarrow A$ and $g :W_q\rightarrow N$, where the brackets $[-,-]_{\underset{K}{\smile}}$ and $[-,-]_{\underset{K}{\frown}}$ denote the graded commutators. In fact, we have
\begin{eqnarray*}\theta_N ([f,g]_{\underset{K}{\smile}}) = \theta_N (f\underset{K}{\smile} g- (-1)^{pq} g\underset{K}{\smile} f)= (-1)^{np}f \underset{K}{\frown} \theta_N (g)-(-1)^{pq} \theta_N (g) \underset{K}{\frown} f \\
=(-1)^{np}(f \underset{K}{\frown} \theta_N (g)-(-1)^{p(n-q)} \theta_N (g) \underset{K}{\frown} f)= (-1)^{np} [f,\theta_N (g)]_{\underset{K}{\frown}}.
\end{eqnarray*}
Now apply (\ref{thetaoncommutator}) when $f=\me_A$ is the fundamental Koszul 1-cocycle and use the fundamental formulas of Koszul calculus (Proposition \ref{fundamental}).

Combining $\theta_N ([\me_A,g]_{\underset{K}{\smile}}) = (-1)^n [\me_A,\theta_N (g)]_{\underset{K}{\frown}}$ with $b_K=- [\me_A, -]_{\underset{K}{\smile}}$ and $b^K = -[\me_A, -]_{\underset{K}{\frown}}$,
we deduce that
\begin{equation} \label{thetaondifferential}
\theta_N (b_K(g))= (-1)^n b^K(\theta_N(g))
\end{equation} 
so that $\theta_N$ is a $(-n)$-degree morphism of complexes, thus a $(-n)$-degree complex morphism of dg $\tilde{A}$-bimodules. Next we replace $N$ by $M$.

It remains to prove that $\theta_M$ is an isomorphism. We shall give an inverse map $\eta_M : M\otimes W_{n-\bullet}\rightarrow  \Hom(W_{\bullet},M)$ to $\theta_M$ in terms of the basis $\{x_1, \ldots, x_n\}$ of $V$. Using the basis (\ref{W_pbasis}) of $W_p$, remark that $\theta_M: \Hom(W_p,M) \rightarrow M\otimes W_{n-p}$ is completely defined by all the images $\theta_M(f)$ where $f:W_p \rightarrow M$ is given by an element $m\in M$ and a list of indices $1 \leq j_1 < \cdots < j_p \leq n$ such that, for all the lists of indices $1 \leq k_1 < \cdots < k_p \leq n$, $f(x_{k_1} \wedge \ldots  \wedge x_{k_p})$ is zero except if $k_1= j_1, \ldots, k_p=j_p$ for which $f(x_{j_1} \wedge \ldots  \wedge x_{j_p})=m$. Then, for such an $f$, Formula (\ref{rightcap1}) provides
$$\theta_M(f)= (-1)^{np} \mathrm{sgn}(\sigma)\, m \otimes (x_{i_1} \wedge \ldots  \wedge x_{i_{n-p}})$$
where the list of indices $1 \leq i_1 < \cdots < i_{n-p} \leq n$ is the complement of $1 \leq j_1 < \cdots < j_p \leq n$ in $\{1, \ldots , n\}$ and $\sigma$ is the corresponding $(p,n-p)$-shuffle. Therefore, it suffices to define $\eta_M : M\otimes W_{n-p}\rightarrow  \Hom(W_p,M)$ by
$$\eta_M (m \otimes (x_{i_1} \wedge \ldots  \wedge x_{i_{n-p}}))=(-1)^{np} \mathrm{sgn}(\sigma)\, f.$$
The last assertion of the theorem is an easy consequence of the specialization $M=A$. 
\Edm

\Bcr \label{dualityiso2}
Keeping notations of Theorem \ref{dualityiso},
$$H(\theta_M)=c \underset{K}{\frown} - : \HK^{\bullet}(A,M) \longrightarrow \HK_{n-\bullet}(A,M)$$
is a $(-n)$-degree isomorphism of graded $\HK^{\bullet}(A)$-bimodules and $H(H(\theta_M)) : \HK_{hi}^{\bullet}(A,M) \rightarrow \HK^{hi}_{n-\bullet}(A,M)$ is a $(-n)$-degree isomorphism of graded $\HK_{hi}^{\bullet}(A)$-bimodules.
\Ecr
\Bdm
The first assertion follows from Theorem \ref{dualityiso}. For the second assertion, we refer to~\cite[Subsection 2.5]{bt:kcpa} for the definition of higher (co)homologies. The $(-n)$-degree isomorphism $H(\theta_M)$ of graded $\HK^{\bullet}(A)$-bimodules satisfies
$$H(\theta_M)(\overline{\me}_A \underset{K}{\smile} \alpha) = (-1)^n \overline{\me}_A \underset{K}{\frown} H(\theta_M)(\alpha)$$
for all $\alpha \in \HK^{\bullet}(A,M)$, where $\overline{\me}_A$ is the class of $\me_A$ in $\HK^1(A)$. The differential of the higher cohomology (homology) is defined on $\HK^{\bullet}(A,M)$ (resp. $\HK_{\bullet}(A,M)$) by $\overline{\me}_A \underset{K}{\smile} - $ (resp. $\overline{\me}_A \underset{K}{\frown} - $). Therefore $H(\theta_M)$ is a $(-n)$-degree isomorphism of complexes for higher (co)homologies. Applying higher (co)homologies, 
$$H(H(\theta_M)): \HK^{\bullet}_{hi}(A,M) \rightarrow \HK^{hi}_{n-\bullet}(A,M)$$
is an $(-n)$-degree graded $\HK^{\bullet}_{hi}(A)$-bimodule isomorphism.
\Edm
\\

The higher Koszul cohomology $\HK_{hi}^{\bullet}(A)$ and the higher Koszul homology $\HK^{hi}_{\bullet}(A)$ were computed in~\cite[Theorem 3.17, Theorem 5.4]{bls:kocal}. The result is the following.
\begin{eqnarray*}
\HK_{hi}^n(A) \cong k \ \ \mathrm{and} \ \ \HK_{hi}^p(A) \cong 0 \ \ \mathrm{if} \ \ p\neq n, \\
\HK^{hi}_0(A) \cong k \ \ \mathrm{and} \ \ \HK^{hi}_p(A) \cong 0 \ \ \mathrm{if} \ \ p\neq 0.
\end{eqnarray*}
The isomorphim $H(H(\theta_A))$ is trivial in this case.

\subsection{Deriving the strong $n$-Kc-Calabi-Yau property for $A=S(V)$} \label{strong property}

\Bte \label{S(V)isstrong}
Let $V$ be a $k$-vector space of finite dimension $n$ and $A=S(V)$ be the symmetric algebra of $V$. Then the quadratic algebra $A$ is strong $n$-Kc-Calabi-Yau.
\Ete

The aim of this subsection is to prove this theorem. Our proof is based on the duality isomorphism in Theorem \ref{dualityiso}.

\Bpo \label{isodgabimod}
Let us choose in Theorem \ref{dualityiso} the $A$-bimodule $M=A^e$ when $A^e$ is naturally considered as a left $A^e$-module. Then the $(-n)$-degree complex isomorphism
$$\theta_{A^e}: \Hom(W_{\bullet},{A^e}) \longrightarrow A^e\otimes W_{n-\bullet}$$
in $\mathcal{C}(\tilde{A})$ can be defined as a $(-n)$-degree complex isomorphism in $\mathcal{C}(\tilde{A}, \abimod)$.
\Epo
\Bdm
As in Part (ii) of Theorem \ref{newobjects}, define $\Hom(W_{\bullet},{A^e})=\Hom_{A^e}(K(A), A^e)$ as an object in $\mathcal{C}(\tilde{A}, \abimod)$. Viewing $A^e$ in $A^e\otimes W_{\bullet}$ as a right $A^e$-module, $A^e\otimes W_{\bullet}$ becomes also an object in $\mathcal{C}(\tilde{A}, \abimod)$. Note that the $A$-bimodule actions $\alpha (m \otimes w) \beta$ on an element $m \otimes w$ in $A^e\otimes W_{\bullet}$ are defined by $(m . (\beta \otimes \alpha)) \otimes w$ for any $\alpha$ and $\beta$ in $A$, where the dot denotes the product in $A^e$.

To conclude that $\theta_{A^e}$ is a morphism in $\mathcal{C}(\tilde{A}, \abimod)$, it suffices to prove that it is a morphism of $A$-bimodules.

As at the end of proof of Theorem \ref{dualityiso}, fix $f: W_p \rightarrow A^e$ by an element $m\in A^e$ and a list of indices $1 \leq j_1 < \cdots < j_p \leq n$ such that, for all the lists of indices $1 \leq k_1 < \cdots < k_p \leq n$, $f(x_{k_1} \wedge \ldots  \wedge x_{k_p})$ is zero except if $k_1= j_1, \ldots, k_p=j_p$ for which $f(x_{j_1} \wedge \ldots  \wedge x_{j_p})=m$. According to Equality (\ref{rightA^elinear}), for any $\alpha$ and $\beta$ in $A$, $\alpha f \beta :W_p \rightarrow A^e$ is defined by the same data as $f$ except that $m$ is replaced by $m.(\beta \otimes \alpha)$. Therefore, one has
$$\theta_{A^e}(\alpha\, f\, \beta)= (-1)^{np} \mathrm{sgn}(\sigma)\, m.(\beta \otimes \alpha) \otimes (x_{i_1} \wedge \ldots  \wedge x_{i_{n-p}})= \alpha \,\theta_{A^e} (f)\, \beta. \ \Edm$$ 

In the previous proof, $A^e\otimes W_{\bullet}$ is defined as an object in $\mathcal{C}(\tilde{A}, \abimod)$. Moreover (see \ref{K(A)dgbimod}) the Koszul complex $K(A)$ is also an object in $\mathcal{C}(\tilde{A}, \abimod)$. Define~\cite[Subsection 3.3]{bt:kcpa} the $0$-degree graded linear map
$$\tilde{\varphi} : A^e \otimes W_{\bullet} \longrightarrow K(A)$$
by $\tilde{\varphi} ((\alpha \otimes \beta) \otimes y_1  \ldots y_q) = \beta \otimes (y_1  \ldots y_q) \otimes \alpha$ for any $\alpha$ and $\beta$ in $A$, where $y_1  \ldots y_q$ is the conventional notation for an arbitrary element of $W_q$. Then the following comes from~\cite[Subsection 3.3]{bt:kcpa}. It is easy to prove it directly.

\Bpo \label{naturaliso}
With the above assumptions, $\tilde{\varphi} : A^e \otimes W_{\bullet} \rightarrow K(A)$ is a $0$-degree complex isomorphism in $\mathcal{C}(\tilde{A}, \abimod)$, that is, an isomorphism in $\mathcal{C}_{str}(\tilde{A}, \abimod)$.
\Epo

Therefore the composite $\tilde{\varphi} \circ \theta_{A^e}: \Hom(W_{\bullet},{A^e}) \longrightarrow K(A)_{n-\bullet}$ is a $(-n)$-degree complex isomorphism in $\mathcal{C}(\tilde{A}, \abimod)$. Shifting the degree in $K(A)$, we see that
\begin{equation} \label{isoA^e}
\tilde{\varphi} \circ \theta_{A^e}: \Hom(W_{\bullet},{A^e}) \longrightarrow K(A)[-n]
\end{equation}
is an isomorphism in $\mathcal{C}_{str}(\tilde{A}, \abimod)$. Consequently,
\begin{equation} \label{HisoA^e}
H(\tilde{\varphi} \circ \theta_{A^e}): \HK^{\bullet}(A,{A^e}) \longrightarrow H(K(A))_{n-\bullet}
\end{equation}
is an isomorphism of $A$-bimodules. Since $A$ is Koszul, $\HK^p(A, {A^e})\cong 0$ if $p\neq n$ and $\HK^n(A, {A^e})\cong A$ as $A$-bimodule. Remark that the differential of $\Hom(W_{\bullet}, {A^e})$ is not zero. Therefore Proposition \ref{symalg} fails for the non central $A$-bimodule $M=A^e$.  

Apply the localization functor $\tilde{Q}:  \mathcal{C}_{str}(\tilde{A},\abimod) \rightarrow \mathcal{D}(\tilde{A},\abimod)$ to the isomorphism (\ref{isoA^e}) and obtain an isomorphism
$$\tilde{Q}(\tilde{\varphi} \circ \theta_{A^e}) : \Hom_{A^e}(K(A), A^e) \longrightarrow K(A)[-n]$$
in $\mathcal{D}^b(\tilde{A},\abimod)$. Conformally to Definition \ref{strongKcCY}, our aim is to prove that there is an isomorphism
$$\RHom_{A^e}(K(A), A^e) \cong K(A)[-n]$$
in $\mathcal{D}^b(\tilde{A},\abimod)$. Thus it suffices to prove that there is an isomorphism
\begin{equation} \label{lastiso}
\Hom_{A^e}(K(A), A^e) \cong \RHom_{A^e}(K(A), A^e)
\end{equation}
in $\mathcal{D}^b(\tilde{A},\abimod)$. The existence of such an isomorphism will be a consequence of some general facts on derived categories applied to derived Koszul calculus. We follow Yekutieli's book~\cite{yeku:dercat} for these general facts.

Actually, we prove in the next proposition the existence of such an isomorphism in the unbounded derived category $\mathcal{D}(\tilde{A},\abimod)$ for \emph{any quadratic algebra $A$ over a finite quiver $\mathcal{Q}$}. In this statement,  there is no finiteness assumption on the length of $K(A)$. If $K(A)$ has a finite length, the isomorphism in $\mathcal{D}(\tilde{A},\abimod)$ induces an isomorphism in $\mathcal{D}^b(\tilde{A},\abimod)$. So Theorem \ref{S(V)isstrong} will be an immediate consequence of the following general result. This result was not included in~\cite{bt:kcpa} but was implicitly used just after~\cite[Definition 5.10]{bt:kcpa} for proving that the preprojective algebras are strong $2$-Kc-Calabi-Yau.

\Bpo \label{generalphi}
Let us go back to the general setup of Subsection \ref{koscalA}. Let $A$ be a quadratic algebra over a finite quiver $\mathcal{Q}$. Then there exists an isomorphism \begin{equation} \label{phi}
\phi: \Hom_{A^e}(K(A), A^e) \longrightarrow \RHom_{A^e}(K(A), A^e)
\end{equation}
in $\mathcal{D}(\tilde{A},\abimod)$.
\Epo
\Bdm
Let us consider the triangulated functor 
$$F:= \Hom_{A^e}(-,A^e) : \mathcal{K}(\tilde{A},\abimod)^{op} \longrightarrow \mathcal{D}(\tilde{A},\abimod)$$
defined in Formula (\ref{threeQtriangfunc3}) of Subsection \ref{basicfunctors}. In this subsection, we have applied a general fact of the theory, namely Theorem \ref{contrarderfunctexist}, since $\mathcal{C}(\tilde{A},\abimod)$ has enough $K$-projectives (Proposition $\ref{fundaresult}$). So $F$ has a triangulated right derived functor
$$(RF, \eta^R):\mathcal{D}(\tilde{A},\abimod)^{op} \longrightarrow \mathcal{D}(\tilde{A},\abimod)$$
where $\eta^R: F \Rightarrow RF \circ Q$ is an isomorphism of triangulated functors from $\mathcal{K}(\tilde{A},\abimod)^{op}$ to $\mathcal{D}(\tilde{A},\abimod)$. Here $Q$ stands for the localization functor from $\mathcal{K}(\tilde{A},\abimod)^{op}$ to $\mathcal{D}(\tilde{A},\abimod)^{op}$. Recall that for every $K$-projective object $P$ in $\mathcal{K}(\tilde{A},\abimod)$, the morphism $\eta^R_P : F(P)  \rightarrow RF\circ Q(P)$ in $\mathcal{D}(\tilde{A},\abimod)$ is an isomorphism.

Define the morphism
$$\phi :=\eta^R_{K(A)} : F(K(A)) \longrightarrow RF \circ Q(K(A))$$
in $\mathcal{D}(\tilde{A},\abimod)$. Note that
$$F(K(A))=\Hom_{A^e}(K(A), A^e) \ \ \mathrm{and} \ \ RF \circ Q(K(A))=\RHom_{A^e}(K(A), A^e).$$
Our aim is to show that $\phi$ is an isomorphism in $\mathcal{D}(\tilde{A},\abimod)$. However we can't conclude that immediately because we don't know whether the object $K(A)$ is $K$-projective in $\mathcal{K}(\tilde{A},\abimod)$ -- may be it is not the case. In fact, the isomorphism $\theta_A$ in Theorem \ref{dualityiso} shows that the free $\tilde{A}$-bimodules seem too big compared to $K(A)$.

Our idea is to work with the underlying $A$-bimodules and to check that $K(A)$ is $K$-projective in $\mathcal{K}(\abimod)$. Let us introduce the dg functor
\begin{equation} \label{functorUnd}
\Und : \mathcal{C}(\tilde{A},\abimod) \longrightarrow \mathcal{C}(\abimod)
\end{equation}
where $\mathcal{C}(\abimod)$ is the dg category of complexes of $A$-bimodules (see Subsection \ref{specializations}). By definition the functor $\Und$ sends any object in $\mathcal{C}(\tilde{A},\abimod)$ described in Subsection \ref{dgbimod} to its underlying $A$-bimodule, and analogously on the morphisms. A general fact~\cite[Theorem 5.5.1]{yeku:dercat} provides a triangulated functor
\begin{equation} \label{functorUndtriangulated}
\Und : \mathcal{K}(\tilde{A},\abimod) \longrightarrow \mathcal{K}(\abimod).
\end{equation}

The commutative diagram
\begin{eqnarray} \label{HUnd1}
 \mathcal{K}(\tilde{A},\abimod) \ \ \ \stackrel{\Und}{\longrightarrow} & \mathcal{K}(\abimod) \nonumber  \\
\downarrow H  \ \ \ \ \ \ \ \ \ \ \   &  \downarrow H \\
\mathcal{G}_{str}(\abimod) \ \ \ \stackrel{id}{\longrightarrow} & \mathcal{G}_{str}(\abimod) \nonumber
\end{eqnarray}
shows that the functor $\Und : \mathcal{K}(\tilde{A},\abimod) \longrightarrow \mathcal{K}(\abimod)$ sends quasi-isomorphisms to quasi-isomorphisms and consequently passes trivially through the derived categories. Since $\Und$ is its own left derived functor, the so-obtained left derived functor is traditionnally still denoted by $\Und$, and one has the commutative diagram
\begin{eqnarray} \label{HUnd2}
 \mathcal{D}(\tilde{A},\abimod) \ \ \ \stackrel{\Und}{\longrightarrow} & \mathcal{D}(\abimod) \nonumber  \\
\downarrow H  \ \ \ \ \ \ \ \ \ \ \   &  \downarrow H \\
\mathcal{G}_{str}(\abimod) \ \ \ \stackrel{id}{\longrightarrow} & \mathcal{G}_{str}(\abimod) \nonumber
\end{eqnarray}
Similarly we have the right derived functor $\Und : \mathcal{D}(\tilde{A},\abimod)^{op} \longrightarrow \mathcal{D}(\abimod)$ satisfying an analogous commutative diagram.

Let us define now the triangulated functor 
$$G:= \Hom_{A^e}(-,A^e) : \mathcal{K}(\abimod)^{op} \longrightarrow \mathcal{D}(\abimod).$$
One has the commutative digagram
\begin{eqnarray} \label{FG}
 \mathcal{K}(\tilde{A},\abimod)^{op} \ \ \ \stackrel{F}{\longrightarrow} & \mathcal{D}(\tilde{A},\abimod) \nonumber  \\
\downarrow \Und  \ \ \ \ \ \ \ \ \ \ \   &  \downarrow \Und \\
\mathcal{K}(\abimod)^{op} \ \ \ \stackrel{G}{\longrightarrow} & \mathcal{D}(\abimod) \nonumber
\end{eqnarray}
 According to Subsection \ref{specializations}, $\mathcal{C}(\abimod)$ has enough $K$-projectives and Theorem \ref{contrarderfunctexist} applies. Thus $G$ has a triangulated right derived functor $RG$ and one has the commutative diagram
\begin{eqnarray} \label{RFRG}
 \mathcal{D}(\tilde{A},\abimod)^{op} \ \ \ \stackrel{RF}{\longrightarrow} & \mathcal{D}(\tilde{A},\abimod) \nonumber  \\
\downarrow \Und  \ \ \ \ \ \ \ \ \ \ \   &  \downarrow \Und \\
\mathcal{D}(\abimod)^{op} \ \ \ \stackrel{RG}{\longrightarrow} & \mathcal{D}(\abimod) \nonumber
\end{eqnarray}
Recall that $\Und$ is identified to its left and right derived functor.

Consider the chain complex $K(A)$ as a cochain complex by setting $K(A)^p=K(A)_{-p}$ for any $p\in \mathbb{Z}$. This cochain complex is bounded above and it is a cochain complex of projective $A$-bimodules. Then $K(A)$ is a semi-projective complex in $\mathcal{C}(\abimod)$~\cite[Proposition 11.3.5]{yeku:dercat}, so that $K(A)$ is $K$-projective in $\mathcal{C}(\abimod)$~\cite[Theorem 11.3.2]{yeku:dercat} and in $\mathcal{K}(\abimod)$ as well. Therefore the morphism in $\mathcal{D}(\abimod)$ defined by
$$\Und(\phi): G(K(A)) \longrightarrow RG \circ Q(K(A))$$
is an isomorphism in $\mathcal{D}(\abimod)$. Here $Q$ is the localization functor from $\mathcal{K}(\abimod)^{op}$ to $\mathcal{D}(\abimod)^{op}$.

Consequently, the commutative diagram (\ref{HUnd2})
shows that $H(\phi)=H(\Und(\phi))$ is an isomorphism in $\mathcal{G}_{str}(\abimod)$. A fundamental fact says to us that the functor
$$H: \mathcal{D}(\tilde{A},\abimod) \longrightarrow \mathcal{G}_{str}(\abimod)$$
is \emph{conservative}~\cite[Corollary 7.2.13]{yeku:dercat}. Therefore, since $H(\phi)$ is an isomorphism, $\phi$ is an isomorphism.
\Edm


\begin{thebibliography}{99}

\bibitem{bls:kocal} R. Berger, T. Lambre, A. Solotar, Koszul calculus, \emph{Glasg. Math. J.} \textbf{60} (2018), 361-399.
\bibitem{rbnm:kogo} R. Berger, N. Marconnet, Koszul and Gorenstein properties for homogeneous algebras, 
  \emph{Algebras and Representation Theory} \textbf{1} (2006), 67-97.
\bibitem{bs:cupcap} R. Berger, A. Solotar, A cup-cap duality in Koszul calculus, arXiv:2007.00627v2.
\bibitem{bt:kcpa} R. Berger, R. Taillefer, Koszul calculus of preprojective algebras, \emph{J. London. Math. Soc.} \textbf{102} (2020), 1241-1292.
\bibitem{bou:alg} N. Bourbaki, \emph{\'{E}l\'ements de math\'ematique, Alg\`ebre, Chapitres 1 à 3}, Hermann, 1970.
\bibitem{ce:homo} H. Cartan, S. Eilenberg, \emph{Homological algebra}, Princeton Mathematical Series \textbf{19}, Princeton Univ. Press, 1956.
\bibitem{dgt:variant} Y. Daletski, I. Gelfand, B. Tsygan, On a variant of noncommutative differential geometry, \emph{Soviet Math. Dokl.} \textbf{40} (1990), 422–426.
\bibitem{eu:calculus} C-H. Eu, The calculus structure of the Hochschild homology/cohomology of preprojective algebras of Dynkin quivers, \emph{J. Pure Appl. Algebra} {214} (2010), 28-46. 
\bibitem{eusched:cyfrob} C-H. Eu, T. Schedler, Calabi-Yau Frobenius algebras, \emph{J. Algebra} {321} (2009), 774-815.
\bibitem{gm:methods} S. Gelfand, Y. Manin , \emph{Methods of homological algebra}, Springer Monographs in Mathematics, Springer-Verlag, 2003.
\bibitem{vg:cy} V. Ginzburg, Calabi-Yau algebras, \texttt{arXiv:math.AG/0612139v3}.
\bibitem{hat:algtop} A. Hatcher, \emph{Algebraic topology}, Cambridge University Press, Cambridge, 2002. 
\bibitem{keller:derdgcat} B. Keller, Deriving DG categories, \emph{Ann. Sci. École Norm. Sup.} \textbf{27} (1994), 63–102.
\bibitem{keller:clustdercat} B. Keller, Cluster algebras and derived categories. Derived categories in algebraic geometry, 123–183, EMS Ser. Congr. Rep., Eur. Math. Soc., Z\"{u}rich, 2012.
\bibitem{tl:bvcy} T. Lambre, Dualit\'e de Van den Bergh et structure de Batalin-Vilkoviski\v{\i} sur les alg\`ebres de Calabi-Yau, \emph{J. Noncommut. Geom.} \textbf{4} (2010), 441-457.
\bibitem{lpv:poisson} C. Laurent-Gengoux, A. Pichereau, P. Vanhaecke, \emph{Poisson structures}, Grundlehren Math. Wiss., \textbf{347} Springer, 2013.
\bibitem{spalt:unbounded} N. Spaltenstein, Resolutions of unbounded complexes,
  \emph{Compositio Math.} \textbf{65} (1988), 121-154.
\bibitem{tt:calculus} D. Tamarkin and B. Tsygan, The ring of differential operators on forms in noncommutative calculus, \emph{Proc. Sympos. Pure Math.} \textbf{73}, Amer. Math. Soc. (2005), 105-131.
\bibitem{toen:dag} B. Toën, \emph{Derived algebraic geometry}, EMS Surv. Math. Sci. (2014), 153–240.
\bibitem{vdb:nch} M. Van den Bergh, Noncommutative homology of some three-dimensional quantum spaces, \emph{K-Theory} {8} (1994), 213-230.
\bibitem{vdb:existence} M. Van den Bergh, Existence theorems for dualizing complexes over non-commutative graded and filtered rings. \emph{J. Algebra} \textbf{195} (1997), 662–679. 
\bibitem{vdb:dual} M. Van den Bergh, A relation between Hochschild homology and cohomology for Gorenstein rings, \emph{Proc. Amer. Math. Soc.} {126} (1998), 1345-1348; erratum \emph{ibid.} {130} (2002), 2809-2810.
\bibitem{weib:homo} C. A. Weibel, \emph{An introduction to homological algebra}, Cambridge Studies in Advanced Mathematics \textbf{38}, Cambridge Univ. Press, 1994.
\bibitem{yeku:dercat} A. Yekutieli, \emph{Derived Categories}, Cambridge Studies in Advanced Mathematics \textbf{183}, Cambridge Univ. Press, 2020.

\end{thebibliography}
\end{document}